\documentclass[a4paper,12pt]{amsart} 
\usepackage[utf8]{inputenc}
\usepackage[T1]{fontenc}
\usepackage[english]{babel}
\usepackage{lmodern}
\usepackage{transparent}
\usepackage{enumitem}
\usepackage{amsmath,amssymb,amscd,amsfonts,latexsym}
\usepackage{hyperref}
\usepackage{caption}
\usepackage{tikz-cd}
\usepackage{subcaption}
\usepackage{wrapfig}
\usepackage{float}
\restylefloat{table}
\usepackage{diagbox}
\usepackage{tikz}
\usepackage{url}
\usepackage{hyperref}
\usepackage{enumitem}
\usepackage{csquotes}
\usepackage{verbatim}
\usepackage{float}
\usepackage{cleveref}
\usepackage{breqn}
\newtheorem{theorem}{Theorem}[section]
\newtheorem{prop}{Proposition}[section]
\newtheorem{lemma}[theorem]{Lemma}

\newtheorem*{theorem*}{Theorem}
\theoremstyle{definition}
\newtheorem{definition}{Definition}[section]

\newtheorem{problem}[theorem]{Problem}
\newtheorem{remark}[theorem]{Remark}

\newtheorem{assumption}[theorem]{Assumption}
\numberwithin{equation}{section}

\date{}

\def \Q {{\mathbb Q}}
\def \C {{\mathbb C}}
\def \H {{\mathbb H}}

\def \N {{\mathbb N}}
\def \Z {{\mathbb Z}}
\def \R {{\mathbb R}}

\def \C {{\mathbb {C}}}

\begin{document}

\title{Higher genus Angel surfaces}

\author[R. Bardhan]{Rivu Bardhan}
\address{Department of Mathematics, Shiv Nadar University, Dadri 201314, Uttar Pradesh, India}
\email{rb212@snu.edu.in}

\author[I. Biswas]{Indranil Biswas}
\address{Department of Mathematics, Shiv Nadar University, Dadri 201314, Uttar Pradesh, India}
\email{indranil.biswas@snu.edu.in, inrdranil29@gmail.com}

\author[S. Fujimori]{Shoichi Fujimori}
\address{Department of Mathematics, Hiroshima University, Higashihiroshima, Hiroshima 739-8526, Japan}
\email{fujimori@hiroshima-u.ac.jp}

\author[P. Kumar]{Pradip Kumar}
\address{Department of Mathematics, Shiv Nadar University, Dadri 201314, Uttar Pradesh, India}
\email{pradip.kumar@snu.edu.in}

\date{}

\subjclass[2020]{53A10, 49Q05, 53C42, 30F60}

\keywords{Minimal surface, finite total curvature, two ends, angel surface}

\begin{abstract}
We prove the existence of complete minimal surfaces in $\mathbb{R}^3$ of arbitrary genus $p\, \ge\, 1$ and least 
total absolute curvature with precisely two ends --- one catenoidal and one Enneper-type --- thereby solving,
affirmatively, a problem posed by Fujimori and Shoda. 
These surfaces, which are called \emph{Angel surfaces}, generalize some
examples numerically constructed earlier by Weber. The construction of these minimal surfaces
involves extending the orthodisk method developed by
Weber and Wolf \cite{weber2002teichmuller}. A central idea in our
construction is the notion of \emph{partial symmetry}, which enables us to introduce controlled
symmetry into the surface.
\end{abstract}

\maketitle

\tableofcontents

\section{Introduction}

Complete minimal surfaces in $\mathbb{R}^3$ with finite total curvature have long been one of the central
objects of 
study in differential geometry. Two topological invariants associated with such surfaces are the 
genus and the number of ends, which are related to the total curvature of the surface
by the \emph{Osserman inequality}~\cite{osserman1964annals}.
It is natural to ask whether there exists a complete minimal
surface, with given genus and number of ends, that minimizes the total absolute curvature.

This question remains open in general, and no complete classification is currently known. In this work, we 
restrict our attention to the simplest nontrivial cases, that is to surfaces with one or two ends. Classical 
examples illustrate the genus-zero case: the plane, with one planar end, and the catenoid, with two 
catenoidal ends. A celebrated result of Schoen, \cite{schoen1983uniqueness}, establishes that the catenoid is 
the unique complete minimal surface of genus zero with two ends and least total absolute curvature.

For higher genus surfaces with one end, Weber and Wolf, \cite{weber1998minimal}, and
independently Sato, \cite{sato1996tohoku}, constructed examples 
realizing the least total absolute curvature. The two-ended case, however, remained largely open. Motivated 
by the known examples and the structure of the Weierstrass representation, Weber provided numerical evidence 
that for each $p \,\ge \, 1$ there exists a complete minimal surface of genus $p$ 
with the least total absolute curvature and exactly two ends --- one asymptotic to an end of
the catenoid and the other to the end of Enneper surface \cite[Example 4.1]{fujimoriangel2016}. Such
surfaces are referred to as \emph{Angel surfaces}.

Progress toward this problem was made by Fujimori and Shoda~\cite{fujimoriangel2016}, who constructed 
explicit examples of genus-one, with a catenoidal end and an Enneper end, and also of even genus, with two twice-wrapped 
catenoidal ends, surfaces which attain the least total absolute curvature. Their construction crucially
relies on the presence of symmetry in order to reduce the complexity of the period problem and also
to control the 
global geometry. However, to date, no such construction is known for any odd genus greater than one.

\begin{figure}[htbp]
  \centering
  \begin{minipage}[b]{0.3\hsize}
    \centering
    \includegraphics[width=\linewidth]{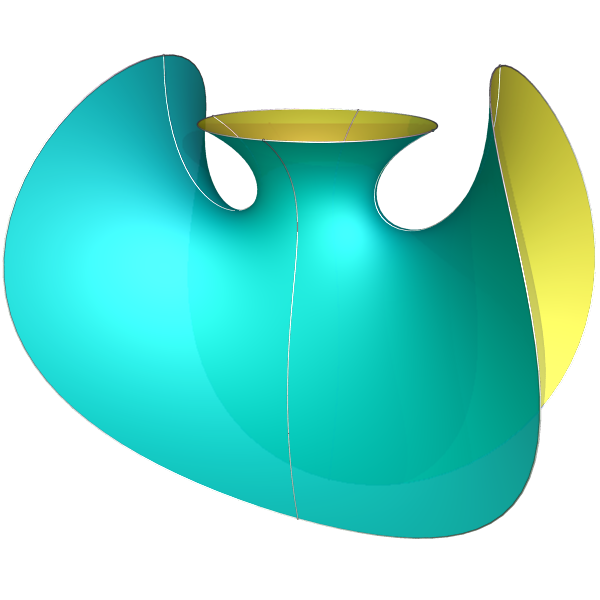}
    \subcaption{$p=0$}\label{angel0}
  \end{minipage} \quad
  \begin{minipage}[b]{0.3\hsize}
    \centering
    \includegraphics[width=\linewidth]{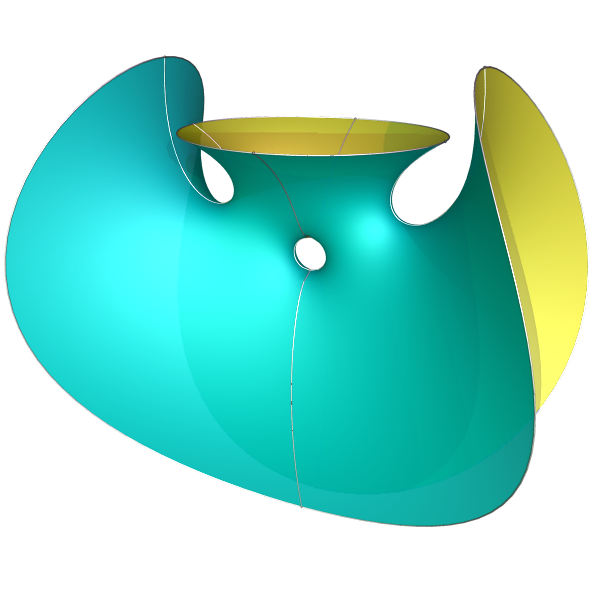}
    \subcaption{$p=1$}\label{angel1}
  \end{minipage} \quad
  \begin{minipage}[b]{0.3\hsize}
    \centering
    \includegraphics[width=\linewidth]{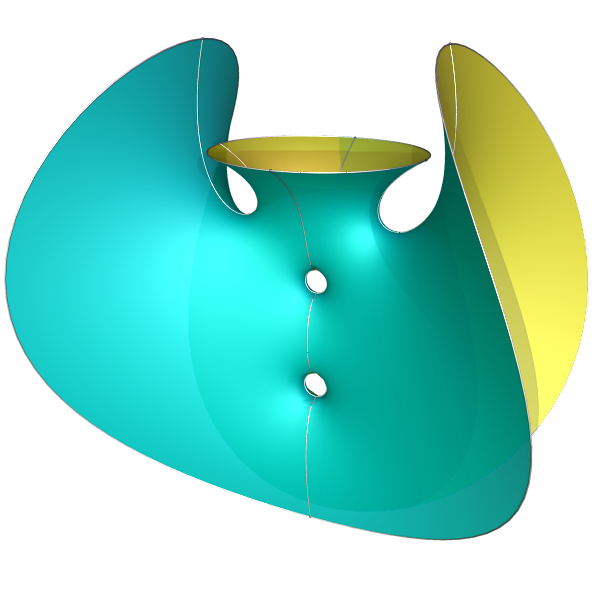}
    \subcaption{$p=2$}\label{angel2}
  \end{minipage} \\
  \begin{minipage}[b]{0.3\hsize}
    \centering
    \includegraphics[width=\linewidth]{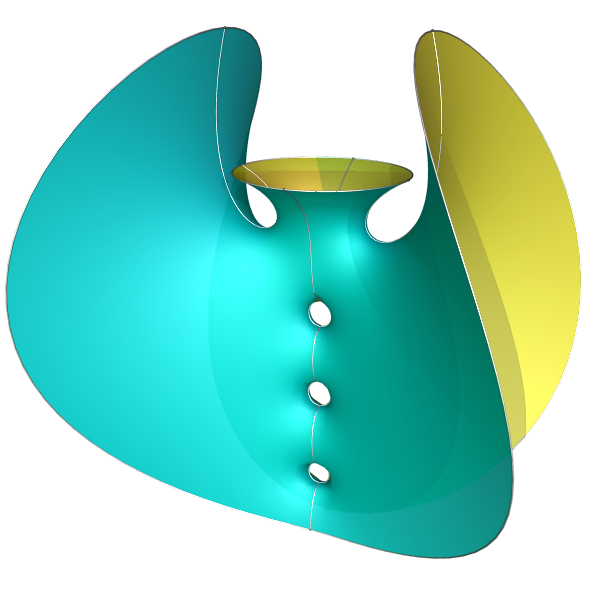}
    \subcaption{$p=3$}\label{angel3}
  \end{minipage} \quad
  \begin{minipage}[b]{0.3\hsize}
    \centering
    \includegraphics[width=\linewidth]{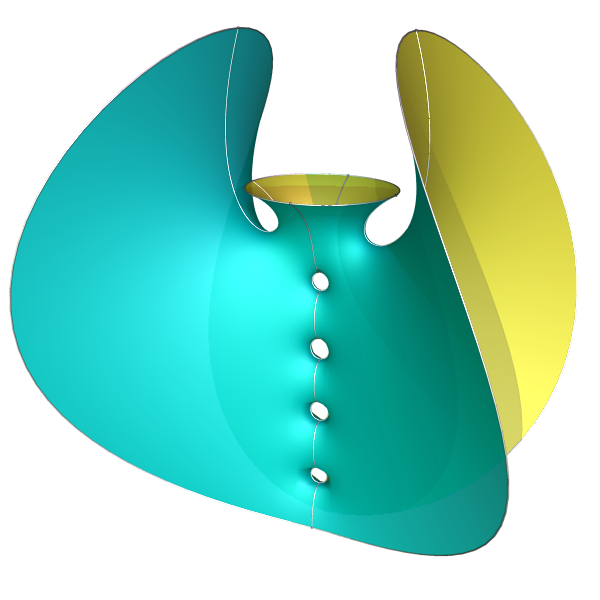}
    \subcaption{$p=4$}\label{angel4}
  \end{minipage}
  \caption{The Angel Surfaces}
  \label{fig:angel}
\end{figure}

The purpose of this article is to construct these surfaces for all genera $p\,\ge\, 1$, thus 
affirmatively solving the problem (Problem~\ref{conjecture:Weber}) in general. In particular, we 
prove that for each $p\,\ge\, 1$ there exists a complete minimal immersion of genus $p$ with least
absolute curvature with one
Enneper end and one catenoid end; see Figure~\ref{fig:angel}. 

For genus one, in~\cite{fujimoriangel2016} Fujimori and Shoda showed using explicit Weierstrass data that 
there exist constants $c\,>\,0$ and $t\,>\,1$ such that on the elliptic curve
\[
\left\{ (z, w) \in (\C\cup\{\infty\})^2 \;\middle|\; w^2 = \frac{z(z-1)}{z - t} \right\}\setminus\Bigl\{(0,0),\,(\infty,\,\infty)\Bigr\}
\]
the pair
$$
G(z,\,w) \ =\ c\,\frac{w}{\,z+1\,}, \ \ \,\, \eta \ =\,\frac{\,z+1\,}{\,z\,}\,dz
$$
constitutes the Weierstrass data of a complete minimal surface of genus $1$ with one Enneper end at 
$(\infty,\, \infty)$ and one catenoid end at $(0,\, 0)$ (see Theorem~\ref{theorem:existence of one genus 
Angel surface by fujimori}). They solved the period problem using symmetry.

To construct Angel surfaces for arbitrary genus $p > 1$, we will ``add handles” to the
above genus 1 Angel surface.
 
The main difficulty lies in solving the period problem for a proposed Weierstrass data set on a higher-genus Riemann surface. More precisely, one seeks a compact Riemann surface $\overline{M}$ of genus $p$, together with a meromorphic function $G$ and a holomorphic 1-form $\eta$ on $M = \overline{M} \setminus \{p_1, p_2\}$, such that

\begin{enumerate}[label=(\roman*)]
  \item The zeros of $\eta$ coincide with the zeros and poles of $G$, i.e.  
  \[
    (\eta)_0 = (G)_0 + (G)_\infty,
  \]
  \item The period condition as in \eqref{equation:period_condition_minimal_surface} is satisfied, and  
  \item At $p_1$ and $p_2$, both $G$ and $\eta$ extend meromorphically, with $p_1$ corresponding to a catenoidal end and $p_2$ to an Enneper end.  
\end{enumerate}

For this, we generalize an approach developed by Weber and Wolf~\cite{weber1998minimal,weber2002teichmuller}
that translates the period problem to a problem in Teichmüller theory. In their pioneering
work~\cite{weber2002teichmuller}, Weber and Wolf introduced the notion of an \emph{orthodisk} and showed
how to encode minimal surface data via Schwarz--Christoffel mappings of planar polygonal domains.
By analyzing moduli of these polygonal domains and using extremal length techniques in Teichmüller
space, they solved the period problem for various families of minimal surfaces. The key idea is to
parametrize the Weierstrass data by certain polygonal domains (orthodisks) in the complex plane, so
that the complicated period integrals on a surface of genus $p$ are converted into conditions on the
geometry of a polygon. By tuning the polygon (and varying it in its moduli), it is possible to achieve
the required period balance.

\begin{figure}[h]
  \centering
  \begin{minipage}{0.40\linewidth}
    \centering
    \def\svgwidth{\linewidth}
\begingroup%
  \makeatletter%
  \providecommand\color[2][]{%
    \errmessage{(Inkscape) Color is used for the text in Inkscape, but the package 'color.sty' is not loaded}%
    \renewcommand\color[2][]{}%
  }%
  \providecommand\transparent[1]{%
    \errmessage{(Inkscape) Transparency is used (non-zero) for the text in Inkscape, but the package 'transparent.sty' is not loaded}%
    \renewcommand\transparent[1]{}%
  }%
  \providecommand\rotatebox[2]{#2}%
  \newcommand*\fsize{\dimexpr\f@size pt\relax}%
  \newcommand*\lineheight[1]{\fontsize{\fsize}{#1\fsize}\selectfont}%
  \ifx\svgwidth\undefined%
    \setlength{\unitlength}{394.04553223bp}%
    \ifx\svgscale\undefined%
      \relax%
    \else%
      \setlength{\unitlength}{\unitlength * \real{\svgscale}}%
    \fi%
  \else%
    \setlength{\unitlength}{\svgwidth}%
  \fi%
  \global\let\svgwidth\undefined%
  \global\let\svgscale\undefined%
  \makeatother%
  \begin{picture}(1,1.05825334)%
    \lineheight{1}%
    \setlength\tabcolsep{0pt}%
    \put(0,0){\includegraphics[width=\unitlength,page=1]{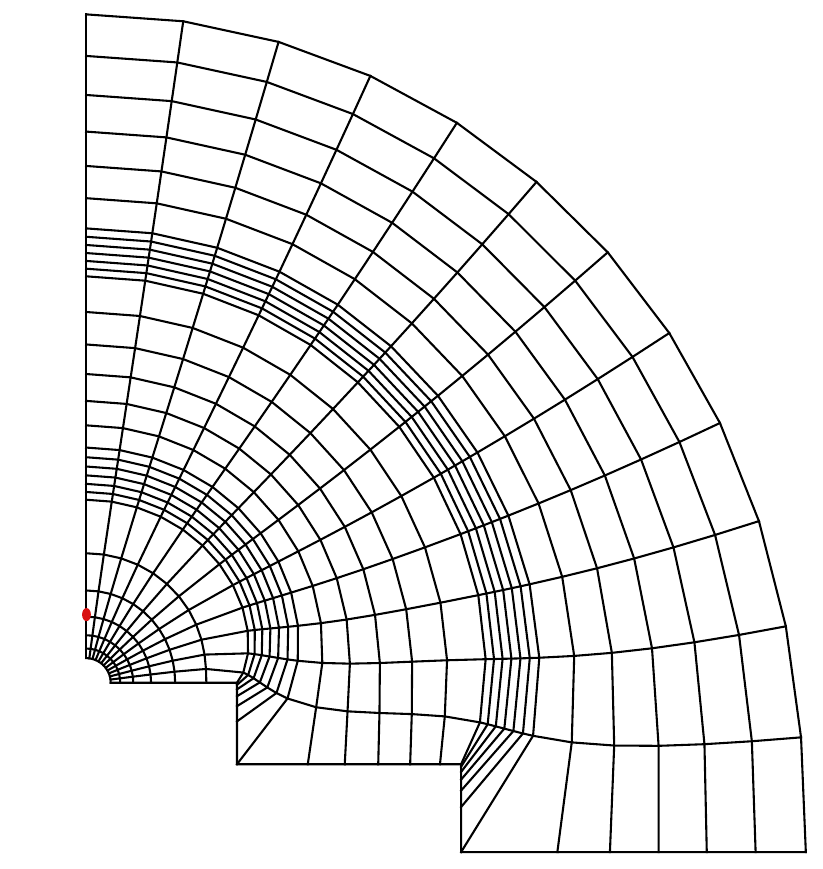}}%
    \put(-0.00167657,0.30026106){\color[rgb]{0.03529412,0,0}\makebox(0,0)[lt]{\lineheight{0}\smash{\begin{tabular}[t]{l}$t_{-1}$\end{tabular}}}}%
    \put(0.21184239,0.19350159){\color[rgb]{0.03529412,0,0}\makebox(0,0)[lt]{\lineheight{0}\smash{\begin{tabular}[t]{l}$t_1$\end{tabular}}}}%
    \put(0.23986675,0.08407311){\color[rgb]{0.03529412,0,0}\makebox(0,0)[lt]{\lineheight{0}\smash{\begin{tabular}[t]{l}$t_2$\end{tabular}}}}%
    \put(0.48141008,0.08407317){\color[rgb]{0.03529412,0,0}\makebox(0,0)[lt]{\lineheight{0}\smash{\begin{tabular}[t]{l}$t_3$\end{tabular}}}}%
    \put(0.47740662,0.01467939){\color[rgb]{0.03529412,0,0}\makebox(0,0)[lt]{\lineheight{0}\smash{\begin{tabular}[t]{l}$t_4$\end{tabular}}}}%
    \put(0,0){\includegraphics[width=\unitlength,page=2]{angel2gn.pdf}}%
  \end{picture}%
\endgroup%

  \end{minipage}
  \hfill
  \begin{minipage}{0.50\linewidth}
    \centering
    \def\svgwidth{\linewidth}
\begingroup%
  \makeatletter%
  \providecommand\color[2][]{%
    \errmessage{(Inkscape) Color is used for the text in Inkscape, but the package 'color.sty' is not loaded}%
    \renewcommand\color[2][]{}%
  }%
  \providecommand\transparent[1]{%
    \errmessage{(Inkscape) Transparency is used (non-zero) for the text in Inkscape, but the package 'transparent.sty' is not loaded}%
    \renewcommand\transparent[1]{}%
  }%
  \providecommand\rotatebox[2]{#2}%
  \newcommand*\fsize{\dimexpr\f@size pt\relax}%
  \newcommand*\lineheight[1]{\fontsize{\fsize}{#1\fsize}\selectfont}%
  \ifx\svgwidth\undefined%
    \setlength{\unitlength}{360bp}%
    \ifx\svgscale\undefined%
      \relax%
    \else%
      \setlength{\unitlength}{\unitlength * \real{\svgscale}}%
    \fi%
  \else%
    \setlength{\unitlength}{\svgwidth}%
  \fi%
  \global\let\svgwidth\undefined%
  \global\let\svgscale\undefined%
  \makeatother%
  \begin{picture}(1,0.81388887)%
    \lineheight{1}%
    \setlength\tabcolsep{0pt}%
    \put(0,0){\includegraphics[width=\unitlength,page=1]{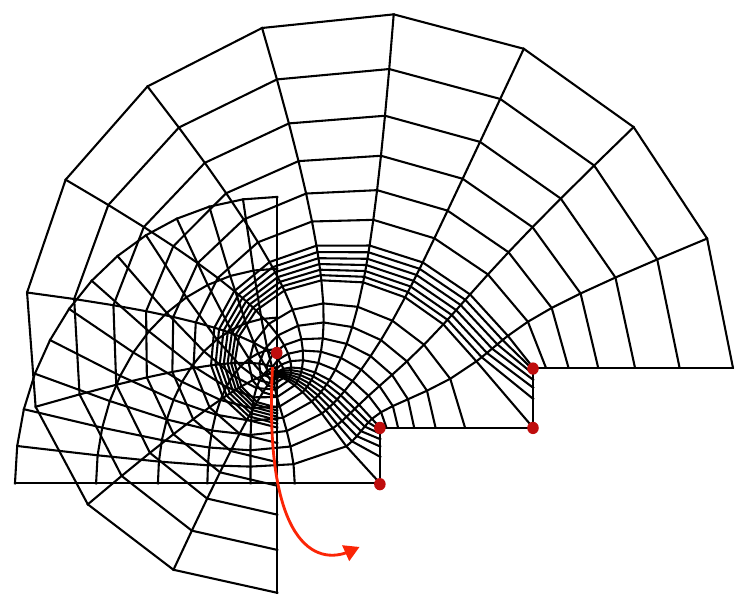}}%
    \put(0.47827517,0.06671219){\color[rgb]{0.03529412,0,0}\makebox(0,0)[lt]{\lineheight{0}\smash{\begin{tabular}[t]{l}$t_{-1}$\end{tabular}}}}%
    \put(0.46801176,0.13445079){\color[rgb]{0.03529412,0,0}\makebox(0,0)[lt]{\lineheight{0}\smash{\begin{tabular}[t]{l}$t_{1}$\end{tabular}}}}%
    \put(0.51624982,0.21142635){\color[rgb]{0.03529412,0,0}\makebox(0,0)[lt]{\lineheight{0}\smash{\begin{tabular}[t]{l}$t_2$\end{tabular}}}}%
    \put(0.67738546,0.20526829){\color[rgb]{0.03529412,0,0}\makebox(0,0)[lt]{\lineheight{0}\smash{\begin{tabular}[t]{l}$t_3$\end{tabular}}}}%
    \put(0.71741269,0.29661274){\color[rgb]{0.03529412,0,0}\makebox(0,0)[lt]{\lineheight{0}\smash{\begin{tabular}[t]{l}$t_4$\end{tabular}}}}%
  \end{picture}%
\endgroup%

  \end{minipage}
  \caption{Genus $2$: flat structures---left for $G\eta$ and right for $G^{-1}\eta$.}
  \label{figure:flat_metrics_due_to_angel_surface_of_genus_2}
\end{figure}

In \cite{weber2002teichmuller}, an orthodisk was defined with help of Schwarz–Christoffel map with odd-integer vertex data so that consecutive boundary edges meet orthogonally and parallel edges come in alternating horizontal/vertical families; the associated flat structures for $G\,\eta$ and $G^{-1}\eta$ are then arranged as a reflexive pair: they share the same conformal polygon (same ordered vertex set) and have conjugate period vectors. This symmetry of the polygonal combinatorics was crucial in their applications (e.g., Costa towers $DH_{m,n}$).  

\begin{figure}[t]
  \centering
  \includegraphics[scale=0.65]{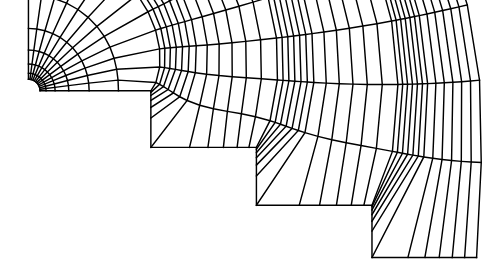}\qquad
  \includegraphics[scale=0.65]{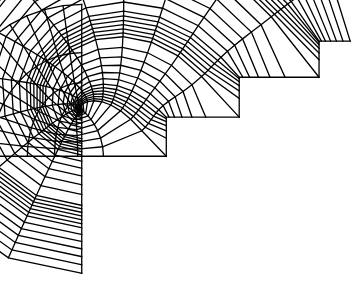}
  \caption{Genus $3$: flat structures---left for $G\eta$ and right for $G^{-1}\eta$.}
  \label{fig:Genus3:flat}
\end{figure}

In contrast, in our setting,  Figure \ref{figure:flat_metrics_due_to_angel_surface_of_genus_2} shows that the two flat (for genus 2) patches determined by $G\,\eta$ and $G^{-1}\eta$ do not share that full combinatorial symmetry: Angel surfaces carry mixed end types (one catenoidal, one Enneper), only partial symmetry (see the staircase in  Figures \ref{figure:flat_metrics_due_to_angel_surface_of_genus_2} and \ref{fig:Genus3:flat}), and hence non-symmetric angle/edge data in the developed polygons. To encode this, we replace the classical notion by essential orthodisk $(c,T,A)$ together with an enhanced conformal polygon $(T_0\subset T)$.  We explained the  modified setup in Section \ref{sec:Generalized orthodisk and minimal surface}.

\subsection{Overview of the construction}

We now give an overview of the main steps.

\subsubsection*{Step 1 -- From orthodisk data to a minimal surface:}\,
We begin by defining a generalized orthodisk $X\, =\, (c,\, T,\, A)$ (Definition~\ref{defn:Generalizedorthodisk})
as a triple, where $$T\ :=\ \{t_1,\, \ldots,\, t_n\,\,\big\vert\,\, t_1 \,<\, t_2\, <\, \cdots 
\,<\, t_n\}\ \subset\ \R,$$ $A$ denotes the tuple $(a_i)_{i=1}^n \,\in \,\mathbb{Q}^n$, and $c\, >
\, 0$ is a constant. Such data determine a Schwarz-Christoffel mapping $F\,:\,\H \,\longrightarrow\, \C$
(see \eqref{equation: Schwarz Christoffel map}) that sends the upper half-plane $\H$ conformally onto a
polygon in the plane with vertices $F(T)$ and interior angles $\pi a_i$ at $F(t_i)$ for each
$t_i \,\in\, T$. The pair $(\H,\,F^*(dz))$ is called a generalized orthodisk.

Let $T_0 \,\subset \, T$ with the cardinality of $T_0$ being odd. Now $(c, \,T,\, T_0,\, A)$ is referred
to as an 
enhanced conformal polygon (see Definition~\ref{defn:essentialConformalPolygon}). As in 
\cite{weber2002teichmuller}, Take the double of $(\H,\,F^*(dz))$ and form a hyperelliptic Riemann surface 
$R_X^{\mathrm{ess}}$ by taking a double cover of the Riemann sphere branched exactly over $T_0 \cup 
\{\infty\}$. By construction, $R_X^{\mathrm{ess}}$ has genus $\frac{\#T_0-1}{2}$. Next, lift the 
differential $F^*(dz)$ to this cover, thus obtaining a meromorphic 1-form $\omega_X$ on $R_X^{\mathrm{ess}}$ 
(Remark~\ref{rem:omegaX}). The divisor of $\omega_X$ on $R_X^{\mathrm{ess}}$ is determined directly by the 
$A$ and $T$.

Now, if one can find two such 1-forms $\omega_{X_1}$ and $\omega_{X_2}$ on the same Riemann surface $R$ that 
are related by complex conjugation of periods, then they can be combined to produce a minimal surface. In 
formal terms, Theorem~\ref{theorem:existence of minimal surface from reflexive orthodisk} loosely states that if 
$(X_1,\, X_2)$ is an e-reflexive pair of genus-$p$ generalized orthodisks 
(Definition~\ref{definition: essentially reflexive orthodisk}), then one can construct a complete 
minimal surface of genus $p$ in $\R^3$.

\subsubsection*{Step 2 -- Constructing genus-$p$ generalized orthodisk data matching the ends:}\, Step 1 shows
that it is needed to find an e-reflexive pair of generalized orthodisks in order to obtain the desired 
minimal surface. As it is evident from the definition, to find such a pair of e-reflexive 
generalized orthodisks, one must begin with a 
pair of e-conjugate generalized orthodisks, as described in Definition~\ref{defn:essentially 
conjugate generalized othodisks}.

We will determine the general form of the required pair $(X_1,\, X_2)$ of genus $p$ by generalizing the pair 
from the known case of genus 1. Intuitively, when constructing the polygons $X_1$ and $X_2$, there is 
flexibility in choosing the lengths and positions of their edges, which in turn affects the period balance. 
We refer to this free parameter as $\lambda$, which essentially controls the “handle size”. As $\lambda$ 
varies, the pair $(X_1,\, X_2)$ moves within a family of e-conjugate genus-$p$ orthodisk pairs. 
Our goal in the final step will be to show that for some value of $\lambda$, the pair becomes reflexive.

\subsubsection*{Step 3 -- Finding a reflexive pair:}\,
The final step is to prove that as the parameter $\lambda$ varies, a reflexive pair indeed arises. To this 
end, we introduce a moduli space of e-conjugate orthodisk pairs with a certain partial symmetry. 
Similar to \cite{weber2002teichmuller}, define on this moduli space a real-valued height function 
$H_p(\lambda)$ (Section~\ref{sec:height_func}) which measures the non-reflexivity of ${X_1}$ and ${X_2}$.

 Following the same strategy as in \cite{weber2002teichmuller}, it is shown that $H_p(\lambda) \ge 0$ for all $\lambda$, and that $H_p(\lambda) = 0$ if and only if $(X_1,\, X_2)$ is a reflexive pair of genus $p$. It is then proved that such a $\lambda = (X_1,\, X_2)$ exists.

\subsection{Organization of the article}

The paper is organized as follows. Section~\ref{sec:Preminaries} reviews preliminaries on minimal 
surfaces and recalls Weber and Wolf’s orthodisk method in its original form (Subsection~\ref{subsection: 
Weber and Wolf's method}). Section~\ref{sec:Generalized orthodisk and minimal surface} introduces the
generalized orthodisk framework. In Section~\ref{sec:Weber's conjecture on higher genus Angel surface},
we state our 
problem (Problem~\ref{conjecture:Weber}), and Section~\ref{sec:Orthodisk for the genus 1case} 
revisits the genus 1 Angel surface and it is reformulated in terms of an orthodisk 
(Proposition~\ref{proposition:genus 1 Angel surface}).
Section~\ref{section:formalWeierstrassdataForGenusP} presents the proposed Weierstrass data for the genus 
$p$ Angel surface. In the remaining sections, we carry out the general methods proposed by Weber and Wolf in 
\cite{weber2002teichmuller}, adapted to our setting.

\medskip

\noindent \textbf{Acknowledgment.}
Figures~1--3 were generated in \textit{Mathematica} with the mesh package developed by Matthias~Weber. We are grateful to Michael~Wolf for his encouragement and feedback; his suggestions improved the exposition.

\section{Preliminaries}\label{sec:Preminaries}

\subsection{Minimal surface in $\mathbb R^3$}
We start by recalling the Weierstrass-Enneper representation for the minimal surfaces in $\mathbb R^3.$

\subsubsection{Weierstrass-Enneper representation}

For an oriented conformal minimal immersion $X\,:\, M\,\longrightarrow\, {\mathbb R}^3$, the induced metric on $M$ 
produces a complex structure on $X$, making it a Riemann surface. This Riemann surface is equipped with a 
meromorphic function $G$ and a holomorphic one-form $\eta$ such that the divisor $\{\eta\,=\, 0\}$ coincides 
(with multiplicities) with the divisor $\{G\,=\, 0\} + \{G\,=\, \infty\}$. Moreover,
\begin{equation}\label{equation 2.1:weirstrass enneper representation for minimal surfaces}
X(p)\ =\ X(x_0)+ \text{Re} \int^p_{x_0} \left(\frac{1}{2} (G^{-1} - G),\, \frac{\sqrt{-1}}{2} ( G^{-1} + G ),\, 1 \right)\eta.
\end{equation}

The above triple $\left( M,\,G,\,\eta\right)$ is referred to as the Weierstrass data for the minimal
surface $(X,\, M)$. Conversely, any such triplet $\left( M,\,G,\,\eta\right)$ for which
$(\eta)_0 \,=\, (G)_0 + (G)_\infty$ and
$$
\mathrm{Re} \int_\gamma \left(\frac{1}{2} (G^{-1} - G),\, \frac{\sqrt{-1}}{2} (G^{-1} + G),\, 1 \right) \eta
\ =\ 0
$$
for all $\gamma \,\in\, H_1(M,\, \mathbb{Z})$, the function in
\eqref{equation 2.1:weirstrass enneper representation for minimal surfaces} defines a conformal minimal immersion.

For an oriented complete minimal immersion with finite total curvature, the corresponding Riemann surface is the
complement of finitely many points of a compact Riemann surface.

\subsubsection{Construction of minimal surfaces in $\R^3$}\label{subsub:construction_Of_minimalSurface}

Take a closed Riemann surface \(\overline{M}\),
and $M\, \subset\, \overline{M}$ with $0\,<\,\# (\overline{M}\setminus M)\, <\, \infty$,
together with a pair \((G,\, \eta)\), where \(G\) is a meromorphic function on \(\overline{M}\) and \(\eta\) is a
holomorphic 1-form on \(M\) which are meromorphic on $\overline{M}$. On $M$, these are required to satisfy the following two conditions:
\medskip

\noindent{(i)  The Divisor Condition.}
$(\eta)_0\ =\ (G)_0 + (G)_\infty$.

\noindent {(ii) The Period Condition.}
This condition is to ensure that the map in \eqref{equation 2.1:weirstrass enneper representation for minimal surfaces}
is well defined. The period condition consists of the following requirements: For every closed curve
\(\gamma \,\in\, H_1(M,\,\mathbb Z)\),
\begin{align}\label{equation:period_condition_minimal_surface}
    \int_{\gamma} G\eta &\,=\, \overline{\int_{\gamma} G^{-1}\eta}, \quad\operatorname{Re} \int_{\gamma} \eta
\,=\, 0.
\end{align}

\subsubsection{The metric and ends of a minimal surface}

Let $\left(M,\,G,\,\eta\right)$ be the Weierstrass data of a minimal surface $X:M\rightarrow\R^3$. Then the Riemannian metric on $M$ induced by the immersion
$X$ in \eqref{equation 2.1:weirstrass enneper representation for minimal surfaces} is
\begin{equation}\label{eds}
ds^2\ =\ \tfrac{1}{4}\left(\lvert G\rvert+\frac{1}{\lvert G\rvert}\right)^2\lvert \eta\rvert^2.
\end{equation}
The completeness of the minimal surface is equivalent to the completeness of $M$ as a Riemannian manifold with
respect to the metric $ds^2$ in \eqref{eds}. At the punctures (their neighborhoods will be called the ends of the minimal
surface), this completeness condition translates into the condition 
\[ds^2\ \sim\ \operatorname{o}\left(\lvert z\rvert^{-k}\right)\text{ for some }k\,\geq\, 4.\]

There are several types of ends based on their asymptotic behavior. Two such ends are as in Table \ref{table 1:type_of_ends}
\begin{table}[ht]
    \centering
    \[
    \begin{array}{| c |c|c|}
        \hline
       \operatorname{Ord}_{(\cdot)} G & \operatorname{Ord}_{(\cdot)} \eta& \text{Type of end}   \\ \hline
         \pm\,1 & -1 &\text{Catenoid end} \\ \hline
        \pm\,1 & -3 & \text{Enneper end} \\ \hline
    \end{array}
    \]
    \caption{Criteria for Catenoid and Enneper ends}
    \label{table 1:type_of_ends}
\end{table}

\subsection{The method of Weber and Wolf to construct minimal surfaces}\label{subsection: Weber and Wolf's method}

We recall from \cite{weber2002teichmuller} the main tools, namely the Schwarz-Christoffel mapping and
the orthodisk.

\subsubsection{Schwarz-Christoffel mapping}\label{pairCTA}

Given a real constant $c\,>\,0,$  an ordered subset $T\,=\,\lbrace t_1,\,t_2,\,\ldots,\,t_n\rbrace
\, \subset\,  \R$ with \(t_1\,<\,t_2\,<\,\ldots\,<\,t_n\) and tuple $A
\,=\,(a_1,\,a_2,\, \ldots ,\, a_n)\in \mathbb Q^n$, a Schwarz-Christoffel mapping is the following:
\begin{equation}\label{equation: Schwarz Christoffel map}
F\,:\, \mathbb{H}\cup\R\cup\lbrace\infty\rbrace\,\longrightarrow\, \C\cup\lbrace\infty\rbrace,\ \ \,
z\, \mapsto\, c\int_{\sqrt{-1}}^{z} \prod_{i=1}^{n} \left( t - t_i \right)^{a_i-1} dt.
\end{equation}
This map $F$ sends the interior of $\mathbb{H}$ biholomorphically to the interior of an Euclidean polygon
--- with possibly infinity as a vertex --- and carries the real line continuously to the boundary of the
polygon, where each $t_i$ is mapped to some vertex of the polygon. The interior angle at $F(t_i)$ is $a_i\pi.$
If $\infty$ is a vertex, define
\[a_\infty\ =\ n-1 \ - \sum_{i=1}^{n}a_i.\]
The interior angle at ~$F(\infty)$~ is $a_{\infty}\pi.$

\subsubsection{Schwarz-Christoffel mapping and orthodisk \cite{weber2002teichmuller}}\label{defn:orthodisk}

Consider the Schwarz-Christoffel mapping $F$ (see \eqref{equation: Schwarz Christoffel map})
corresponding to an ordered subset $T\,=\, \lbrace t_1,\,t_2,\,\dots,\,t_n\rbrace$ of $\R,$
~$A\,=\, \left(a_i\right)_{i=1}^n$, $2a_i\,\in\,2\Z+1$,
and $c\,=\,1.$ Consider $\H\cup\R$, and equip $\H$ with the pullback $F^*(dx^2+dy^2)$
of the Euclidean metric. This pair $(\H\cup\R,\, F^*(dx^2+dy^2))$ is called an ~\textit{orthodisk}.
The elements of $T$ are called the \textit{vertices} of the orthodisk, and $A$ is called the
\textit{vertex data} for the orthodisk. The Riemann surface with boundary, namely, $\H\cup\R\cup\lbrace\infty\rbrace$, together
with the marked points at $T$, is called the \textit{conformal polygon} of the orthodisk.

We observe that an orthodisk is uniquely identified with the triplet $(\H,\,T,\,A).$ Therefore, from now 
onwards, orthodisks are denoted through their corresponding triples $(\H,\,T,\,A).$
An orthodisk will often be identified with its image in $\mathbb{C}$.

The edges of an orthodisk are the boundary
segments between vertices; they come in a natural order. Consecutive edges
meet orthogonally at the finite vertices. Every other edge is parallel (cf. \cite[Definition 3.1.1]{weber2002teichmuller}) for the flat metric of the orthodisk. Oriented distances
between parallel edges are called periods. The periods can have four different
classes: $+1,\, -1,\, +\sqrt{-1},\, -\sqrt{-1}$.

Two orthodisks are said to be conformal if they share the same conformal polygon. So
$(\H,\, T,\, A)$ and $(\H,\, T',\, A')$ are conformal if and only if $T\,=\, T'$.

Two orthodisks with equal number of vertices~say, \( X_1 \) and \( X_2 \) are said to be \textit{conjugate}~ if there is a line in \(\C\) so  that the corresponding periods are 
symmetric with respect to that line.

Finally, reflexive orthodisks are defined as follows: two orthodisks are said to be \textit{reflexive} if  they are both conformal and conjugate.

\subsubsection{Minimal surface from a reflexive pair} \label{subsec:minimalsurfacefromReflexivepair} 
Given an orthodisk $X\,=\, (\H,\, \{t_j\}_{-p}^{p},\,(a_j)_{-p}^{p})$, Weber and Wolf constructed a
hyperelliptic Riemann surface $R_X$ which is branched over the vertices of the
representative of the conformal polygon of the
orthodisk in the double of $\H\cup\R$. Pull back to $R_X$ the unique meromorphic form corresponding to the 
metric of the orthodisk; call this pullback as $\omega_X$.
For an orthodisk $X\,=\,\left( \H,\,\{t_j\},\,(a_j)\right)$, if the vertex data at $t_i$ is $a_i$, then the
corresponding angle at the 
image of $t_i$ --- which is the Euclidean angle at $F(t_i)$ (see \eqref{equation: Schwarz Christoffel map})
--- is $a_i\pi.$ Further the order of $\omega_X$ at the representative of $t_i$ (including infinity)
in $R_X$ is $2a_i-1.$

If $X_1=(\H,\,\{t_j\}_{-p}^{p},\,A=(a_j)_{-p}^p),\, X_2=(\H,\,\{t_j\}_{-p}^p,\,B=(b_j)_{-p}^p)$ are conformal orthodisks, then the definition of conformal orthodisk ensures that although as 
flat surfaces $\left(R_{X_1},\,\omega_{X_1}\right),\, \left(R_{X_2},\,\omega_{X_2}\right)$ might be different,
nevertheless $R_{X_1}$ and $R_{X_2}$ carry the same conformal structure. In other words, the underlying Riemann surfaces coincide.

 We have $R\,=\,R_{X_{1}}\,=\, R_{X_{2}}$ for any pair of reflexive orthodisks $\left(X_{1},\,X_{2}\right)$.
Let $\omega_{X_1}$ and $\omega_{X_2}$ be the corresponding meromorphic forms on $R_{X_1}$ and
$R_{X_2}$ respectively (see previous paragraph).  For notational convenience, denote \(X_1\) by \(X_{G\eta}\), \(X_2\) by \(X_{G^{-1}\eta}\), $\omega_{X_1}$ by $\omega_{G\eta}$ and $\omega_{X_2}$ by $\omega_{G^{-1}\eta}$.  This renaming indicates that
for the obtained minimal surface, $\omega_{G\eta}$ (respectively,
$\omega_{G^{-1}\eta}$) plays the role of $G\eta$ (respectively, $G^{-1}\eta$) on the Riemann surface $R_{X_{G\eta}}$
as mentioned in the period condition (see \eqref{equation:period_condition_minimal_surface}).

Further, if we take 
\[\omega\ =\ \prod{\left(t-t_i\right)^{\frac{a_i+b_i}{2}-1}}dt\ \text{ and }\ \pi_{G\eta}\ :\
R_{X_{G\eta}}\ \longrightarrow\ \C,\] 
where \(a_i\) (respectively \(b_i\)) is the vertex data for the orthodisk \(X_1\) (respectively \(X_2\)) with \(2a_i,\,2b_i\in 2\Z+1\) and \(a_i\,+\,b_i\,\equiv\,0\pmod{2}\), 
then  it is proved in \cite[Theorem 3.3.5]{weber2002teichmuller}  that the triplet  \[\left(R_{X_{G\eta}},\ G\,=\,\frac{\omega_{G\eta}}{\eta},\ \eta
\,=\,\pi_{G\eta}^*(\omega)\right)\] defines a minimal surface when \((X_{G\eta},\,X_{G^{-1}\eta})\) is a reflexive pair.

In \cite{weber2002teichmuller}, the existence of a reflexive pair of orthodisks
for the Costa towers type minimal surfaces is proved.

The following basis of homology of a  hyperelliptic Riemann surface  with infinity as a branch point
will be used in solving the period problem and also to parameterize the polygon space.

\subsubsection{Homology basis of the  hyperelliptic Riemann surface}\label{subsection:homology basis for Riemann 
surface} For $p \,\in\,\N,$ 
$$M^{\text{hyp}}_{p}\ =\ \lbrace\,(z,\,w)\,\,\big\vert\,\, w^2\,=\,
\prod\limits_{j=0}^{2p}(z-t_j)\rbrace$$ is a
 hyperelliptic Riemann surface  of genus $p$ where $\{t_j\}_{j=0}^{2p}$ are real numbers in increasing 
order. Construct a canonical basis $\{A_j,\, B_j\}_{j=0}^{p-1}$
of $H_1(M^{\text{hyp}}_p,\, \mathbb{Z})$ as follows:

\begin{enumerate}
    \item \textbf{$A_j$-cycles (encircling branch cuts):} 
    For each $j\, =\, 0,\, \dots,\, p-1$, the cycle $A_j$ is defined as the lift to the upper
sheet of a small counterclockwise loop in $\mathbb{C}$ that encloses the interval $[t_{2j},\,t_{2j+1}]$
exactly once, while avoiding all other branch points.

    \item \textbf{$B_j$-cycles (connecting adjacent cuts):} 
    For $j \,=\, 0,\, \dots,\, p-1$, the clockwise cycle $B_j$ is constructed as follows:
    \begin{itemize}
        \item Start on the upper sheet at a point just to the left of $t_{2j+1}$ on the real
line, i.e., on the $(j+1)^{\text{th}}$ branch cut,

\item follow a path in $\{z\,\in\,\mathbb{C}\, \big\vert\,\, \operatorname{Im}(z) \,>\, 0\}$ to a point
just to the right of $t_{2j+2}$ on the real line,

\item cross to the lower sheet through the $(j+2)^{\text{th}}$ branch cut,

\item return along the reflection of the initial path in $\{z\,\in\,\mathbb{C}\, \big\vert\,\,
\operatorname{Im}(z) \,<\, 0\}$ of second sheet.
 \end{itemize}
\end{enumerate}

\section{Generalized orthodisk and minimal surface}\label{sec:Generalized orthodisk and minimal surface}
We now present a few definitions and results that generalize the work of Weber and Wolf
recalled in Section \ref{subsection: Weber and Wolf's method}.  
\begin{definition}[{Generalized orthodisk}]\label{defn:Generalizedorthodisk}
A generalized orthodisk is the pair \((\H,\,F^*(dz)),\) where \(F\) is a Schwarz--Christoffel mapping with vertex set \(T\), vertex data \(A,\) and constant \(c>0\) as in Subsection \ref{pairCTA}. For notational convenience, denote a generalized orthodisk by its corresponding tuple \((c,\,T,\,A)\).  
\end{definition}
For the triplet $(c,\,T,\,A)$ (see Subsection \ref{pairCTA}), where $A\,=\, \left({a_i}\right)_{i=1}^n$ for 
$2a_i\,\in\,2\Z+1$ and $c\,=\,1$, the pair $(\mathbb H,\,F^*(dz))$ (see \eqref{equation: Schwarz Christoffel map}
for $F$) is an orthodisk (see Subsection \ref{defn:orthodisk}).

As before, the elements  in $T$ are called the \emph{vertices} of the generalized orthodisk, and
$\mathbb H$ together with $T$ is the \emph{conformal polygon}, while $A$ is the \emph{vertex data}.

\begin{definition}[Enhanced conformal polygon]\label{defn:essentialConformalPolygon}
Let $(c,\,T,\,A)$ be a generalized orthodisk. Take a non-empty ordered subset  $T_0\,\subset\, T$, which
will be referred to as the set of marked vertices.  The set $\mathbb{H} \cup \mathbb{R}\cup\{\infty\}$, together with the marked vertices, is called an \emph{enhanced conformal polygon} (e-conformal polygon) of the generalized orthodisk. For notational convinence sometimes \(T_0\) is referred as e-conformal polygon. Furthermore, $(c,\,T,\,T_0,\,A)$ is referred to as an \emph{enhanced generalized orthodisk} (e-generalized orthodisk).
\end{definition}

Next, we define the notion of conformal e-generalized orthodisk pair.
\begin{definition}[Conformal pair of e-generalized orthodisks]
	Let \(X=(c_1,\,T_{1},\,T_{0,1},\,A_{1})\)
	and \(Y=(c_2,\,T_{2},\,T_{0,2},\,A_{2})\) be e-generalized orthodisks with
	\(
	\#T_{1} \;=\; \#T_{2},\;\#T_{0,1} \;=\; \#T_{0,2}.
	\)
	We call $(X,Y)$ a \emph{conformal pair of e-generalized orthodisks} if there exists a conformal self-map
	\(
	\phi:\mathbb{H}\cup\mathbb{R}\longrightarrow\mathbb{H}\cup\mathbb{R}
	\)
	such that, \(\phi\) takes the ordered set \(T_1\) to  ordered set \(T_2\).
\end{definition}

Notation and definitions are put together in Table \ref{table:generalized-orthodisk-notation} for convenience.

\begin{table}[ht!]
\centering

\begin{tabular}{|p{0.9cm}|p{3cm}|p{11cm}|}
\hline
\textbf{S. No.} & \textbf{Symbol / Term} & \textbf{Description / Meaning} \\
\hline
1 & $(\mathbb{H}, F^*(dz))$ & Generalized orthodisk, defined on upper half-plane $\mathbb{H}$ with pulled back differential. \\
\hline
2 & $(c, T, A)$ & Alternative notation for a generalized orthodisk: $c > 0$ (scale), $T$ (vertices), $A \in \mathbb{Q}^n$ (vertex data). \\
\hline
3 & $T$ & Set of vertices of the conformal polygon. \\
\hline
4 & $\mathbb{H} \cup T$ & Conformal polygon domain. \\
\hline
5 & $A$ & Vertex data; determines angle or order at each $t_j \in T$. \\
\hline
6 & $T_0 \subset T$ & marked subset of vertices used in defining an enhanced conformal polygon. \\
\hline
7 & $(c,\, T,\, T_0,\, A)$ & Enhanced generalized orthodisk; includes marked vertex subset $T_0$. \\
\hline
8 & $\omega = F^*(dz)$ & Meromorphic 1-form on $\mathbb{H}$, extended to $\mathbb{C}$ via reflection. \\
\hline
9 & $R_X^{\mathrm{ess}}$ & Hyperelliptic double cover of $\mathbb{C}$ branched over $T_0 \cup \{\infty\}$. here $\# T_0$ is odd. \\
\hline
10 & $\omega_X = \pi_X^*(\omega)$ & pulled back of $\omega$ to $R_X^{\mathrm{ess}}$. \\
\hline
11 & $P_t$, $P_{t^\pm}$ & Preimages of vertex $t \in T$ under the cover $\pi_X$, depending on whether $t$ is a branch point. \\
\hline
12 & $(X, Y)$ & Pair of e-conjugate generalized orthodisks with conjugate period. \\
\hline
13 & E\,-\,reflexive pair & Two e-generalized orthodisks which are both conformal and e-conjugate. \\
\hline

\end{tabular}
\caption{Notation and Terminology Related to Generalized Orthodisks}
\label{table:generalized-orthodisk-notation}
\end{table}

\begin{remark}\label{rem:omegaX}
Fix an enhanced generalized orthodisk $X\,=\,(c,\,T,\,T_0,\,A)$ with enhanced conformal polygon $T_0$ such that $\#\,T_0$ is odd.  Pull back the Euclidean differential $dz$ by the Schwarz-Christoffel map $F$ to obtain a meromorphic $1$-form on $\mathbb H$, which extends using reflection, to a meromorphic $1$-form (still denoted $\omega$) on the complex plane $\C$.  Now let \[R_X^{\mathrm{ess}}\,=\,\{(z,\,w)\,\in\,(\C\cup\{\infty\})^2\,\,\big\vert\,\,
w^2\,\,=\prod_{t\in T_0}(z-t)\},\] which is a  compact Riemann surface of genus $\frac{\#\,T_0-1}{2}$.
Consider the projection
$$\pi_X\ :\ R_X^{\mathrm{ess}}\ \longrightarrow\ \C,\ \ \,  (z,\,w)\,\mapsto\, z,$$
and denote
$$
\omega_X \ :=\ \pi_X^*(\omega).
$$
Thus $\omega_X$ is a meromorphic $1$-form on $R_X^{\mathrm{ess}}$ whose zeros and poles lie above the 
vertices in  $T\cup\{\infty\}$. In particular, every point $t\,\in\, T$ has either one or two preimages in 
$R_X^{\mathrm{ess}}$ (one if $t\,\in\, T_0$ is a branch point and two otherwise); at each such preimage, one can 
compute the cone angle of $\omega_X$ and the order of zero (or pole) of $\omega_X$ in terms of the 
vertex data $A$.
\end{remark}

\begin{remark}[On notation for homology cycles]\label{rem:homology_basis_notation}
In what follows, the same notation $\{A_j,\, B_j\}_{j=0}^{p-1}$ (as in
Subsection~\ref{subsection:homology basis for Riemann 
	surface}) is used in denoting a canonical homology basis on different hyperelliptic Riemann surfaces of genus $p$. This slight abuse of notation is justified since, in each context, the underlying hyperelliptic surface is determined unambiguously by the set of Weierstrass branch points (cf. Remark~\ref{rem:omegaX}) that are placed \textbf{in an order}. In particular, when the branch points are specified with order and these are the same in number, the corresponding canonical cycles $A_j$ and $B_j$ are well-defined, and no confusion arises.
\end{remark}

\begin{definition}[E-conjugate generalized orthodisks of genus 
$p$]\label{defn:essentially conjugate generalized othodisks}
Consider two e-generalized orthodisks $X_1 \,=\, (c_1,\, T,\, T^p_0,\, A)$ and $X_2 \,=\, (c_2,\, S,\, S^p_0,
\, B)$. Suppose $\#\,T \,=\, \#\,S$ and $\#\,T^p_0 \,=\, \#\,S^p_0 \,=\, 2p + 1$.  $(X_1,\, X_2)$ is said to
be a pair of \emph{e-conjugate generalized orthodisks} of genus \(p\) if
\begin{equation}\label{equation:e-conjugate_condition}
	\int_{A_j} \omega_{X_1}
	= \overline{\int_{A_j} \omega_{X_2}},
	\qquad
	\int_{B_j} \omega_{X_1}
	= \overline{\int_{B_j} \omega_{X_2}},
	\quad \text{for }j\,=\,0,\,1,\,\ldots,\,p.
\end{equation}
\end{definition}
Finally, we define \textit{e-reflexive generalized orthodisks of genus \(p\).}
\begin{definition}[E\,-\,reflexive generalized orthodisks of genus $p$]\label{definition: essentially reflexive orthodisk}
A pair of e-conjugate generalized orthodisks of genus \(p\), say \((X_1,\,X_2)\) is called e-reflexive if it is also a conformal pair.
\end{definition}
\begin{remark}\label{remark: regarding cone angle and order of zero of generalized orthodisk}
Let $X\,=\, (c,\, T,\, T_0,\, A)$ be an enhanced generalized orthodisk with enhanced polygon $T_0$. Suppose that
$t_0 \,\in \,T \setminus T_0$ is a non-marked vertex with corresponding vertex data $a_0$. Under the
Schwarz-Christoffel map, the interior angle at $t_0$ is $a_0\,\pi$. Since $t_0$ is not a branch point of the
hyperelliptic surface $R_X^{\mathrm{ess}}$, it lifts to two distinct points --- which are denoted by
$P_{t_0^+}$ and $P_{t_0^-}$  --- on $R_X^{\mathrm{ess}}$. At each of these points, the differential $\omega_X$ has a cone angle of $2a_0\,\pi$ and a zero (or pole) of order $a_0 - 1$.

In contrast, if $t \,\in\, T_0$ is an marked (branch) vertex with vertex data $a_t$, then it lifts to a single 
point $P_t \,\in\, R_X^{\mathrm{ess}}$. At this point, $\omega_X$ has a cone angle of $4\,a_t\,\pi$ and a 
zero (or pole) of order $2a_t - 1$.

Thus, the cone angles and the divisor of the differential $\omega_X$ on the surface $R_X^{\mathrm{ess}}$
are entirely determined by the vertex data $A$.
\end{remark}

\subsection{From e-reflexive orthodisks to minimal surfaces}
Let $(X_1^p,\,X_2^p)$ be an e-reflexive pair of generalized orthodisks of genus $p$. 
Then there exist two ordered sets of real numbers with equal cardinality and increasing order,
\(
T_1=\{t_1,t_2,\ldots,t_{2p+k}\}\), \(T_2=\{s_1,s_2,\ldots,s_{2p+k}\},
\) 
for some $k\in\mathbb{N},$ 
and there exists a conformal map
\(
\phi:\mathbb{H}\cup\mathbb{R}\longrightarrow\mathbb{H}\cup\mathbb{R}
\)
satisfying
\[
\phi(t_i)=s_i,\;\; \text{for all } i\,=\,1,\,\ldots,\,2p+k.
\]
Furthermore, there exist
\(T^{p}_{0,1}=\{t_{i_1},t_{i_2},\ldots,t_{i_{2p+1}}\}\subseteq T_1\) and 
 \(T^{p}_{0,2}=\{s_{i_1},s_{i_2},\ldots,s_{i_{2p+1}}\}\subseteq T_2,
\)
with \(i_1\,<\,i_2\,<\,\cdots\,<i_{2p+1}.\) Write the orthodisks
\(
X^p_1\ =\ (c_1,\, T_1,\, T^p_{0,1},\, A),\,
X^p_2\ =\ (c_2,\, T_2,\, T^p_{0,2},\, B)
\)
for some \(A,\,B\,\subset\Q^{2p+k}\) and \(c_1,\,c_2\,\in\,\R^{+}.\)
For such data the corresponding hyperelliptic covers of genus $p$ are conformal to each other and there exists a vertex-preserving conformal map
$$
\Phi:\;R^{\mathrm{ess}}_{X^p_1}\;\longrightarrow\;R^{\mathrm{ess}}_{X^p_2}
$$
where $R^{\mathrm{ess}}_{X_j}$ is the surface of genus $p$ branched over $T^p_{0,j} \cup \{\infty\}$ for \(j=1,\,2\). Assume that the vertex data $A\,
=\, (a_1,\, \ldots,\, a_{2p+k})$ and $B\, =\, (b_1,\, \ldots,\, b_{2p+k})$ satisfy the condition $a_j\,+\,b_j\,\equiv\, 0 \pmod{2}
$ for all $j$. Then, on the double of $\mathbb{H} \cup \mathbb{R}$ with marked points
$T_1 \,=\, \{t_1,\, \ldots,\, t_{2p+k}\}$, the meromorphic differential is defined by
$$
\zeta \,= \,\pm \sqrt{c_1 c_2} \;\prod_{i=1}^{2p+k} (t - t_i)^{\frac{a_i\,+\,b_i}{2}-1} \, dt.
$$
The pullback of $\zeta$ to the surface $R^{\mathrm{ess}}_{X^p_1}$, namely $
\eta_{X^p_1}\, =\, \pi_{X^p_1}^*(\zeta),
$  has purely imaginary periods, because all residues are real. Further, since $\Phi$ is the extension of $\phi$ and $\phi$ maps the real line to the real line and preserves the order, it follows that $\Phi$ maps the preimages of the real line in $R^{\mathrm{ess}}_{X^p_1}$ to the preimages of the real line in $R^{\mathrm{ess}}_{X^p_2}$. 
Since the quadratic differentials $\omega_{X^p_1}\,\Phi^*\omega_{X^p_2}$ and $\eta_{X^p_1}^2$ have the same divisor, they differ by a nonzero constant $c_0^2\in\C^\times$, i.e.
\[
\frac{\omega_{X^p_1}\,\Phi^*\omega_{X^p_2}}{\eta_{X^p_1}^2}=c_0^2.
\]
Next, choose $x_0\,\in\,\R$ such that $x_0\,>\, t_{2p+k}$. 
On the preimages of the real line in each $R^{\mathrm{ess}}_{X_j}$, the forms $\omega_{X_j^p}$ and $\eta_{X^p_1}$ take real values for \(j=1,\,2\), and $\Phi$ preserves these preimages. Moreover, since \(\phi\) is obtained via composition of Schwarz-Christoffel maps, \(\phi\vert_{\R}\) is monotonically increasing.
Therefore,
\[
\frac{\omega_{X^p_1}\,\Phi^*\omega_{X^p_2}}{\eta_{X^p_1}^2}(x_0,w(x_0))
\]
is a positive real number for \((x_0,\,w(x_0))\) in \(R^{\mathrm{ess}}_{X^p_1}\) (See Remark \ref{rem:omegaX}). 
Since this value equals the constant $c_0^2$, it follows that $c_0\in\R$.

This leads to the following result:

\begin{lemma}\label{lemma:existence of suitable n}
Take an e-reflexive pair of orthodisks $(X^p_1,\,X^p_2)$
of genus $p$, and let their vertex data be $A$ and $B$ respectively.  If $a_i+b_i\,\equiv\, 0\,\pmod{2}$ for
every $i$, then there exists a meromorphic form $\eta_{X^p_1}$ on $R_{X^p_1}^{\mathrm{ess}}$ and a real constant \(c_0\) such that
\begin{enumerate}
    \item $\displaystyle\operatorname{Re}\int_{\sigma}c_0\eta_{X^p_1}\,=\,0$ for every $\sigma\,\in\,
\pi_1\!\bigl(R_{X^p_1}^{\mathrm{ess}}\bigr)$,

\item $\omega_{X^p_1}\,\Phi^*\omega_{X^p_2}\,=\,c_0^2\eta_{X^p_1}^2$.
\end{enumerate}
\end{lemma}

Note that the pair $(X^p_1,\,X^p_2)$ in Lemma \ref{lemma:existence of suitable n} is e-reflexive; in particular,
this pair is conjugate by definition.  Therefore  the meromorphic forms $\omega_{X^p_1}$ and $\omega_{X^p_2}$ satisfy  the condition \eqref{equation:e-conjugate_condition}, which implies
\[\int_\sigma\omega_{X_1^p}\,=\,\overline{\int_{\sigma}\Phi^*\omega_{X_2^p}}\quad\text{for }\sigma\in\,H_1(R_{X^p_1}^{\mathrm{ess}}).\]
Combining all these, we have the following:

\begin{theorem}\label{theorem:existence of minimal surface from reflexive orthodisk}
	Let $(X_1^p,\,X_2^p)$ be an $e$-reflexive pair of genus-$p$ orthodisks whose corresponding vertex data have even sums.  
	Let $R^{\mathrm{ess}}_{X_j^p}$ be the associated hyperelliptic Riemann surfaces and let $\omega_{X_j^p}$ be the corresponding meromorphic $1$-forms for $j=1,2$, as in Remark~\ref{rem:omegaX}. Let $P_1,\dots,P_r$ be the poles of $\eta_{X_1^p}$ on $R^{\mathrm{ess}}_{X_1^p}$, and let $\Phi\colon R^{\mathrm{ess}}_{X_1^p}\longrightarrow R^{\mathrm{ess}}_{X_2^p}$ be the corresponding vertex-preserving conformal map.  
	Assume that
	\[
	\int_{\sigma} \omega_{X_1^p}
	=
	\overline{\int_{\sigma} \Phi^{*}\omega_{X_2^p}}
	\]
	for every loop $\sigma$ encircling an end.
	
	Then there exists $c_0\in\R^{\times}$ such that the Weierstrass data
	\[
	\left(
	R^{\mathrm{ess}}_{X_1^p}\setminus\{P_1,\dots,P_r\},\;
	G=\frac{\omega_{X_1^p}}{c_0\,\eta_{X_1^p}},\;
	\eta=c_0\,\eta_{X_1^p}
	\right)
	\]
	define a minimal surface via the Weierstrass--Enneper representation \eqref{equation 2.1:weirstrass enneper representation for minimal surfaces}.
\end{theorem}

\section{Higher genus angel surface}\label{sec:Weber's conjecture on higher genus Angel surface}

Weber provided numerical evidence that for each genus $p \,\geq\, 1$, there exists a complete minimal
surface of genus $p$ in $\mathbb{R}^3$ with finite total curvature, one Enneper end, and one
catenoidal end \cite{fujimoriangel2016}, \cite{weberwebsite}. These surfaces, which generalize the
genus-one construction by Fujimori and Shoda, are referred to as \emph{Angel surfaces} \cite{weberwebsite}.

\subsection{Riemann surface}

For $p\in\N,$ and \(\,0 \,<\, t_1 \,<\, t_2 \,< \,t_3\, <\, \cdots \,<\, t_{2p} \) define
\[
f_1^p(z)\ =\ \prod_{j=1}^{p}(z - t_{2j-1}), \qquad
f_2^p(z)\ =\ \prod_{j=1}^{p}(z - t_{2j}).
\]
Then the solution set of the equation  \( f(z, w) \,=\, f_2^p(z)w^2 - zf_1^p(z) \)  in \( \C^2 \) is a nonsingular affine algebraic curve, and its
compactification is the compact Riemann surface
\[
M_{p}\ =\ \left\{ (z,\, w) \,\in\, {(\C \cup \{\infty\})}^2 \,
\big\vert\,\, f(z,\, w) \,=\, 0 \right\},
\]
which may equivalently be expressed as
\[
M_{p}\ =\ \left\{ (z,\, w) \,\in\, {(\C \cup \{\infty\})}^2 \,\,\big\vert\,
\, w^2 \,=\,  z\frac{\,f_1^p(z)}{f_2^p(z)} \right\}.
\]

\begin{problem}\label{conjecture:Weber}
Do there exist $c\,>\,0,\,0\,<\,t_1\,<\,t_2\,<\,\cdots\,<\,t_{2p}$ such that the following Weierstrass data constitutes a minimal surface of genus $p$ with one Enneper end and one catenoid end?
    \[
  \left(M_{p}\setminus\left\{(0,\,0),\,(\infty,\,\infty)\right\},\,\, G\,=\,\frac{cw}{z+1},\,
\,\eta\,=\, \frac{z+1}{z}dz\right)
    \]
    where the catenoid end is at $(0,\,0)$ and the Enneper end is at $(\infty,\, \infty).$
\end{problem}

\subsection{The problem in the polynomial model of a hyperelliptic surface}

The hyperelliptic curve
\[
M^{\text{hyp}}_p\ :=\ \left\{(z,\, w) \,\, \big\vert\,\, w^2 \,=\, zf_1^p(z)f_2^p(z)\right\}
\]
is isomorphic to $M_{p}$ by the map
\begin{equation}\label{equation:isomorphism between polynomial and rational model}
    \Psi_p \,\,\colon\,\,  M_p\, \rightarrow\, M^{\text{hyp}}_p, \quad (z,\,w)
\,\mapsto\, (z,\, f_2^p(z)w).
\end{equation}
Using this, Problem~\ref{conjecture:Weber} can be reformulated as follows.

\begin{problem}\label{conjecture:WeberHyp}
Do there exist constants $c \,>\, 0$ and ordered points $0 \,<\, t_1 \,<\, t_2
\,<\, \cdots \,<\, t_{2p}$ such that the following Weierstrass data define a minimal surface of
genus $p$ with one Enneper end and one catenoid end?
\[
\left( M^{\text{hyp}}_p \setminus \{(0,\,0), (\infty,\,\infty)\},\,\, G\circ\Psi_{p}^{-1}
\,=\, \frac{cw}{f_2^p(z)(z+1)}, \,\, (\Psi_{p}^{-1})^*\eta \,=\, \frac{z+1}{z}dz \right),
\]
where:
\begin{itemize}
\item The catenoid end is located at $(0,\,0),$
\item The Enneper end is located at $(\infty,\,\infty).$
\end{itemize}
\end{problem}

\section{Generalized orthodisk for the genus \texorpdfstring{$1$}{1} case}\label{sec:Orthodisk for the genus 1case}

Before proposing the formal data for Angel surfaces of arbitrary genus, it is instructive to revisit the 
genus-one case. This case was previously studied by Fujimori and Shoda \cite{fujimoriangel2016}, 
who provided an explicit construction of a complete minimal surface of genus 1 with one Enneper end and one 
catenoid end.
Their method used carefully chosen Weierstrass data on a genus-one elliptic curve and leveraged 
symmetries to resolve the period problem. Our goal in this section is to reformulate their construction in 
the language of generalized orthodisks introduced earlier, which will then serve as a model for higher-genus 
constructions done in the subsequent sections.
We start with the precise result they established.

\begin{theorem}[Fujimori--Shoda {\cite{fujimoriangel2016}}]\label{theorem:existence of one genus Angel surface by fujimori}
There exist $c\,>\,0$ and $t\,>\,1$ such that on \[
M_1 \ =\ \left\{(z,\, w) \,\in\, {(\C\cup\{\infty\})}^2 \,\,\big\vert\,
\, w^2 \,=\, \frac{z(z - 1)}{z - t} \right\},
\] the meromorphic function $G\,=\,\frac{cw}{z+1}$ and the $1$--form
$\eta\,=\,\frac{z+1}{z}dz$ give a complete  minimal surface of genus $1$  of least absolute curvature with one Enneper end at
$(\infty,\, \infty)$ and one catenoid end at $(0,\,0).$   
\end{theorem}

Consider the map $\Psi_1$ in \eqref{equation:isomorphism between polynomial and rational model}. Define
the elliptic curve $$M^{\text{hyp}}_1\ :=\
\left\{(z,\,w)\,\in\,{(\C\cup\{\infty\})}^2\,\,\,\big\vert\,\,w^2
\,=\,{z(z-1)(z-t)}\right\},$$ and  the meromorphic function $G_0\, =\,G\circ\Psi_1^{-1}(z,\,w)
\,=\, \frac{cw}{(z+1)(z-t)}$ on it as well as the $1$--form
$\eta_0\,=\, (\Psi_1^{-1})^*\eta\,=\, \frac{z+1}{z}dz$. Note that $M_1^{\mathrm{hyp}}$ is an elliptic curve, not a hyperelliptic curve, but we will use this notation consistently for high genus cases.  Theorem
\ref{theorem:existence of one genus Angel surface by fujimori} can be rewritten as follows.

\begin{prop}\label{proposition:genus 1 Angel surface}
There exist $c\,>\,0,\,t\,>\,1$ such that $\left(M_1^{hyp}\setminus \{(0,\,0),\,(\infty,\, \infty)\},\,\, G_0,\,
\,\eta_0\right)$ gives a complete minimal surface of genus $1$ of least absolute curvature with one Enneper end at $(\infty,\,\infty)$ and
one catenoid end at $(0,\,0).$
\end{prop}

The goal is to construct a pair of e-reflexive orthodisks of genus 1 such that the corresponding meromorphic one-forms, as described in Remark \ref{rem:omegaX} and Subsection  ~\ref{subsec:minimalsurfacefromReflexivepair}, are given by $G_0\eta_0$ and $G_0^{-1}\eta_0$. 
To this end, the divisors of these one-forms are computed first. 

\subsection{Divisors of $G_0\eta_0$, $G_0^{-1}\eta_0$, and $\eta_0$}
Define the following distinguished points on the elliptic surface:
$$
P_{-1^{\pm}} = (-1,\, \pm\sqrt{-2(1+t)}), \quad 
P_0 = (0,\, 0), \quad 
P_1 = (1,\, 0), \quad 
P_t = (t,\, 0), \quad 
P_{\infty} = (\infty,\, \infty).
$$
The meromorphic one-forms under consideration are:
$$
G_0\eta_0\, =\, c\,\frac{w}{z(z-t)}\,dz, \qquad 
G_0^{-1}\eta_0 \,=\, \frac{1}{c}\frac{(z+1)^2(z-t)}{wz}\,dz.
$$
Apart from the points $P_{-1^\pm},\, P_0, P_1, P_a,$ and $P_\infty$, these forms have no other zeros or poles. Their divisors are given by:
\begin{equation}\label{equation:divisor_gn_genus1}
(G_0\eta_0)\,\, =\,\, P_{-1^{\pm}}^0\, P_0^0\, P_1^2\, P_t^0\, P_\infty^{-2},
\end{equation}
\begin{equation}\label{equation:divisor_gninv_genus1}
(G_0^{-1}\eta_0)\,\, =\,\, P_{-1^{\pm}}^2\, P_0^{-2}\, P_1^0\, P_t^2\, P_\infty^{-4},
\end{equation}
\begin{equation}\label{equation:divisor_eta0_genus1}
(\eta_0)\,\, =\,\, P_{-1^{\pm}}^1\, P_0^{-1}\, P_1^1\,P_t^1\, P_\infty^{-3}.
\end{equation}

\begin{table}[ht]
\centering
\begin{tabular}{|c|c|c|c|c|c|}
\hline
\textbf{Meromorphic form} & $P_{-1^{\pm}}$ & $P_0$ & $P_1$ & $P_t$ & $P_{\infty}$ \\ \hline
$G_0\eta_0$ & $0$ & $0$ & $2$ & $0$ & $-2$ \\ \hline
$G_0^{-1}\eta_0$ & $2$ & $-2$ & $0$ & $2$ & $-4$ \\  \hline
$\eta_0$ & $1$ & $-1$ & $1$ & $1$ & $-3$ \\\hline
\end{tabular}
\caption{Divisors of $G_0\eta_0$, $G_0^{-1}\eta_0$, and $\eta_0$ for genus $1$}
\label{table:divisors_genus1}
\end{table}

\subsection{Constructing appropriate generalized orthodisk for $G_0\eta_0,\,G^{-1}_0\eta_0$}\label{subsec:genus 1 orthodisk construction}
 From the construction of the Riemann surface $R^{\mathrm{ess}}_X$ and the associated 1-form $\omega_X$, as described in Subsection~\ref{subsec:minimalsurfacefromReflexivepair}, the following data corresponding to the enhanced generalized orthodisks can be extracted. 

Let $t \,>\, 1$ be as in Proposition~\ref{proposition:genus 1 Angel surface}, and define:
\begin{align*}
&T_{G_0 \eta_0} = \lbrace -1, 0, 1, t\rbrace, &
A_{G_0 \eta_0} &= \left( 1, \tfrac{1}{2}, \tfrac{3}{2}, \tfrac{1}{2} \right), &
T^1_{0,G_0 \eta_0} &= \lbrace 0, 1, t \rbrace, \\
&T_{G_0^{-1} \eta_0} = \lbrace -1, 0, 1, t \rbrace, &
A_{G_0^{-1} \eta_0} &= \left( 3, -\tfrac{1}{2}, \tfrac{1}{2}, \tfrac{3}{2} \right), &
T^1_{0,G_0^{-1} \eta_0} &= \lbrace 0, 1, t \rbrace.
\end{align*}
These define the following pair of enhanced generalized orthodisks:
\begin{itemize}
\item $X_{G_0\eta_0}\ =\ \left(c,\, T_{G_0 \eta_0},\, T^1_{0,G_0 \eta_0},\, A_{G_0 \eta_0} \right),$
\item $X_{G_0^{-1}\eta_0}\ =\ \left(\frac{1}{c},\, T_{G_0^{-1} \eta_0},\, T^1_{0, G_0^{-1} \eta_0},\, A_{G_0^{-1} \eta_0} \right).$
\end{itemize}
In what follows, it is verified that this pair of orthodisks indeed corresponds to the Angel surface of
genus 1 introduced in Proposition~\ref{proposition:genus 1 Angel surface}.

 By comparing $T_{G_0\eta_0}$, $T_{G^{-1}_0\eta_0}$ and the marked vertices $T^1_{0,G_0\eta_0}$, $T^1_{0,G^{-1}_0\eta_0}$ of the orthodisks $X_{G_0\eta_0}$ and $X_{G_0^{-1}\eta_0}$, it is deduced that they share the same conformal polygon as well as the same enhanced conformal polygon. Moreover, Proposition
\ref{proposition:genus 1 Angel surface} gives the following: For all
$$
\sigma \,\in\, H_1(M_1^{\text{hyp}};\, \mathbb{Z}) \,=\, H_1(R^{\mathrm{ess}}_{X_{G_0\eta_0}};\, \mathbb{Z}),
$$
$$
\int_{\sigma} G_0\eta_0\ =\ \overline{\int_{\sigma} G^{-1}_0\eta_0}, \ \ \ 
\operatorname{Re} \int_{\sigma} \eta_0\ =\ 0.
$$
This implies that the two generalized orthodisks are e-conjugate. Since their vertex sets and marked vertex sets also coincide, the pair $(X_{G_0\eta_0},\, X_{G_0^{-1}\eta_0})$ forms an  e-reflexive pair. Moreover, since in this case $\phi=\operatorname{Id}$, the constant $c_0$ in Lemma~\ref{lemma:existence of suitable n} is equal to $1$.
 Hence, by Theorem~\ref{theorem:existence of minimal surface from reflexive orthodisk}, the data
$$
\left(R_{X_{G_0\eta_0}} \setminus \lbrace P_0,\, P_{\infty} \rbrace,\ 
G \,=\, \tfrac{\omega_{X_{G_0\eta_0}}}{\eta_{X_{G_0\eta_0}}},\ 
\eta_{X_{G_0\eta_0}}\right)
$$
define a minimal surface.  Furthermore, by comparing the orders of zeros and poles of $G$ and
$\eta$ it is seen that the point $P_0$ corresponds to a catenoid end and $P_{\infty}$ to an Enneper end.

We have the identifications
$$
R_{X_{G_0\eta_0}}\, =\, M^{\mathrm{hyp}}_1, \ \ \,
\omega_{X_{G_0\eta_0}}\, =\, \pm G_0\eta_0, \ \ \,
\omega_{X_{G_0^{-1}\eta_0}}\, =\, \pm G_0^{-1}\eta_0, \ \ \,
\eta_{X_{G_0\eta_0}} = \pm \eta_0.
$$
Therefore, the minimal surface constructed from this orthodisk data is isometric to the genus 1 Angel surface described in Proposition~\ref{proposition:genus 1 Angel surface}.

\subsection{Associated polygonal picture for genus 1 orthodisks}
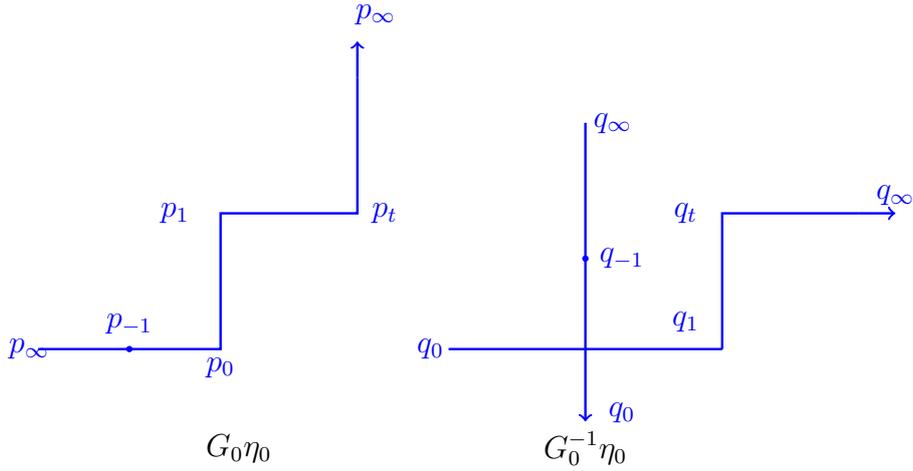
\begin{figure}[ht]
	\centering
	\begin{tikzpicture}[scale=1.2, line width=0.9pt]
		
		\draw[blue] (-2,0) -- (0,0) -- (0,1.5) -- (1.5,1.5) -- (1.5,3);
		\draw[->,blue](1.5,3) -- (1.5,3.4);  
		\fill[blue] (-1,0) circle (1pt) node [above=1pt] {$p_{-1}$};
		\node[blue] at (-2.1,0) {$p_\infty$};
		\node[blue] at (0,-0.2) {$p_0$};
		\node[blue] at (-0.5,1.5) {$p_{1}$};
		\node[blue] at (1.8,1.5) {$p_{t}$};
		\node[blue] at (1.7,3.7) {$p_{\infty}$};
		\node at (0.2,-1.1) {$G_0 \eta_0$};
		
		\draw[blue] (4,0) -- (4,2.5);    
		\node[blue] at (4.3,2.5) {$q_{\infty}$};
		
		\draw[blue] (2.5,0) -- (5.5,0);
		\node[blue] at (2.3,0) {$q_0$};
		
		\draw[blue] (5.5,0) -- (5.5,1.5) -- (7,1.5);
		\draw[->,blue] (7,1.5) -- (7.4,1.5);  
		
		\node[blue] at (5.1,0.3) {$q_{1}$};
		\node[blue] at (5.1,1.5) {$q_{t}$};
		\node[blue] at (7.4,1.7) {$q_{\infty}$};
		\draw[->,blue] (4,0) -- (4,-0.8);
		\fill[blue] (4,1) circle(1pt) node[right=1pt] {\(q_{-1}\)};
		\node[blue] at (4.4,-0.7) {$q_0$};
		
		\node at (4,-1.1) {$G_0^{-1} \eta_0$};
		
	\end{tikzpicture}
	\caption{Image of a pair of generalized orthodisks for genus $1.$ }
	\label{fig:genus 1 orthodisk}
\end{figure}
By the associated polygonal picture of the orthodisk $(c,\, T,\, T_0,\, A)$, we mean the image of $\mathbb{H} \cup \mathbb{R}$ in $\mathbb{C}$ under the map given in Equation~\eqref{equation: Schwarz Christoffel map}.

To draw pictures of the orthodisks $X_{G_0\eta_0}$ and $X_{G_0^{-1}\eta_0}$, the following facts and conventions are used:

\begin{enumerate}
\item The Schwarz-Christoffel map for $X_{G_0\eta_0}$ is given by
$     F_{G_0\eta_0}(z) = c \int_{\sqrt{-1}}^z x^{-\frac{1}{2}}(x - 1)^{\frac{1}{2}}(x - t)^{-\frac{1}{2}}\, dx,
    $
and for $X_{G_0^{-1}\eta_0}$, it is
$     F_{G_0^{-1}\eta_0}(z) = \frac{1}{c}\int_{\sqrt{-1}}^z (x + 1)^2 x^{-\frac{3}{2}}(x - 1)^{-\frac{1}{2}}(x - t)^{\frac{1}{2}}\, dx.
    $
\item Only the boundary is drawn, according to where it maps, with the image of the points $T = \{-1,\, 0,\, 1,\, t\}$ indicated.

\item The convention followed matches the literature, particularly \cite{weber2002teichmuller}. In both orthodisk images, the region $\mathbb{H}$ is mapped to the left side (northwest side) of the boundary image in $\mathbb{C}$.

\item The divisors of \(\omega_{X_{G_0\eta_0}}\) and \(\omega_{X_{G_0^{-1}\eta_0}}\) are obtained by pulling back \(dz\) under the corresponding maps \(F_1\) and \(F_2\). Hence, the divisor data from the previous subsection determines the angles at the images of the points \(t_i \in T\).
\item Note that the illustration of the pair of e-reflexive generalized orthodisks of genus one is first drawn in the complex plane and then rotated so that the resulting polygons have periods (see Remark \ref{rem:omegaX} and Definition \ref{defn:essentially conjugate generalized othodisks}) symmetric with respect to the line \(y \,=\, -x\), in accordance with the convention in the literature. This change in visual representation does not affect the mathematics, as all figures are ultimately mapped to the affine plane \(\C\) in a later section by quotienting with orientation-preserving isometries of \(\R^2\).
\end{enumerate}

Under the corresponding Schwarz-Christoffel maps, label the vertices of the orthodisk corresponding to $G_0\eta_0$ as $p_{-1},\, p_0,\, p_1,\, p_t,$ and $p_\infty$, and those for $G_0^{-1}\eta_0$ as $q_{-1},\, q_0,\, q_1,\, q_t,$ and $q_\infty$ (See Figure \ref{fig:genus 1 orthodisk}).

For $X_{G_0\eta_0}$, all vertices except the one at infinity are mapped to finite points. Consequently, the image of the real line under the Schwarz-Christoffel map is properly embedded in $\mathbb{C}$ without self-intersections.

On the other hand, for $X_{G_0^{-1}\eta_0}$, both $0$ and $\infty$ are mapped to $\infty$. Therefore, the image of the real line under the Schwarz-Christoffel map self-intersects in $\mathbb{C}$. Computing the cone angle at infinity shows that the line connecting the images of $0$ and $\infty$ is parallel to the line segment $q_{1}q_{t}$. Hence, it must intersect either $q_0q_{1}$ or $q_{t}q_\infty$.
\subsection{A genus-one e-reflexive orthodisk pair with unit scaling}\label{subsec:unit-scale}

Recall from Subsection~\ref{subsec:genus 1 orthodisk construction} that the genus--\(1\) Angel surface is
encoded by an e-reflexive pair
\[
\bigl(X_{G_0\eta_0},\,X_{G_0^{-1}\eta_0}\bigr),
\]
whose Schwarz--Christoffel maps (see \eqref{equation: Schwarz Christoffel map} and the explicit formulas in
Subsection~5.3) carry prefactors \(c\) and \(1/c\), respectively.  For the higher genus construction in
Section~\ref{section:formalWeierstrassdataForGenusP}, it is convenient to work with \emph{unit scaling},
i.e.\ representatives for which the scaling constant in \eqref{equation: Schwarz Christoffel map} is \(1\).

Let
\[
(P_1,P_2)
:=\bigl(F_{G_0\eta_0}(\H\cup\R),\,F_{G_0^{-1}\eta_0}(\H\cup\R)\bigr)
\]
be the polygonal images of the pair \(\bigl(X_{G_0\eta_0},X_{G_0^{-1}\eta_0}\bigr)\).

Each \( P_i\) is realized by a Schwarz--Christoffel map of the form
\eqref{equation: Schwarz Christoffel map} with scaling constant \(1\).  Concretely, there exist real
parameters
\[
t_{-1}<t_0<t_1<t_2,
\qquad
s_{-1}<s_0<s_1<s_2,
\]
and Schwarz--Christoffel maps \(\tilde F_{G_0\eta_0},\tilde F_{G_0^{-1}\eta_0}\) with unit scaling, mapping
\(\H\cup\R\) onto \(\tilde P_1\) and \(\tilde P_2\), respectively, such that the marked vertices correspond
to the three finite vertices as in the genus--\(1\) construction.

Define the following unit-scaled enhanced generalized orthodisks (with the same vertex data as in
Subsection~\ref{subsec:genus 1 orthodisk construction}):
\[
\tilde X_{G_0\eta_0}
:=\bigl(1,\;\{t_{-1},t_0,t_1,t_2\},\;\{t_0,t_1,t_2\},\;A_{G_0\eta_0}\bigr),\]
\[
\tilde X_{G_0^{-1}\eta_0}
:=\bigl(1,\;\{s_{-1},s_0,s_1,s_2\},\;\{s_0,s_1,s_2\},\;A_{G_0^{-1}\eta_0}\bigr).
\]

Since the developed polygon for  \(\tilde X_{G_0\eta_0}\) and \(\tilde X_{G_0^{-1}\eta_0}\) are same.  Hence we have the following: 
\begin{lemma}\label{lemma:existence_of_genus_one_orthodisk}
The pair \(\bigl(\tilde X_{G_0\eta_0},\,\tilde X_{G_0^{-1}\eta_0}\bigr)\) is e-reflexive.
\end{lemma}

\section{Proposed data for genus-\texorpdfstring{$p$}{p} Angel surfaces}
\label{section:formalWeierstrassdataForGenusP}
Motivated by the above genus 1 example,  we now present explicit \textbf{formal} Weierstrass data \(\bigl(G,\eta\bigr)\) on a hyperelliptic curve of genus \(p\) which, 
assuming the period conditions hold, will produce a complete minimal surface with one Enneper end and one 
catenoid end. 
\medskip
\subsection{Data for the \(G\eta\)--enhanced generalized orthodisk.}
We choose real parameters
    \(
      t_{-1}\, <\, t_0 \,<\, t_1 \,<\, t_2 \,<\, \cdots \,<\, t_{2p}
    \)
    and set
    \[
      T_{G\eta} \,=\, \{t_{-1},\,t_0,\,t_1,\, \dots,\, t_{2p}\}, \ \ \,
      A_{G\eta}\, = \, \left(1,\, \tfrac12,\, \tfrac32,\,\tfrac12,\,\dots,\, \tfrac32,\, \tfrac12\right),
      \ \ \,
      T_{0,G\eta}^p \,=\, \{t_0,\,t_1,\, \dots,\, t_{2p}\}.
    \]
The e-generalized orthodisk is taken as
    \begin{equation}\label{eqn:generalized orthodisk 1}
      X^p_{G\eta} \ =\  \bigl(1,\;T_{G\eta},\;T_{0,G\eta}^p,\;A_{G\eta}\bigr).
    \end{equation}
    The Schwarz-Christoffel map that maps the upper half plane~ to $\mathbb C$ is given by 
    \begin{equation}\label{eqn:FpG_ETA}
      F_p^{G\eta}(z)
      \ :=\ 
       \int_{\sqrt{-1}}^z
        (t - t_{-1})^{0}
        (t - t_0)^{-\frac12}
        \prod_{k=1}^{2p}(t - t_k)^{(-1)^{k+1}\frac12}
      \,dt.
    \end{equation}

\medskip

\subsection{Data for the \(G^{-1}\eta\)--enhanced generalized orthodisk.} Similarly, the following is taken:
    \(  s_{-1} < s_0 < s_1 < s_2 < \cdots < s_{2p},
    \)
    and set
    \[
      T_{G^{-1}\eta} = \{\,s_{-1},\,s_0,\,s_1,\dots,s_{2p}\}, 
      \quad
      A_{G^{-1}\eta} = \left(\,3,\,-\tfrac12,\,\tfrac12,\,\tfrac32,\,\dots,\, \tfrac12,\, \tfrac32\,\right),
      \quad
      T_{0,G^{-1}\eta}^p = \{\,s_0,\,s_1,\,\dots,\,s_{2p}\}.
    \]
 The e-generalized orthodisk is taken as
    \begin{equation}\label{eqn:genrelized orthodisk2}
      X^p_{G^{-1}\eta} 
      = (1,\;T_{G^{-1}\eta},\;T_{0,G^{-1}\eta}^p,\;A_{G^{-1}\eta}).
    \end{equation}
    The  Schwarz-Christoffel map that maps the upper half plane~ to $\mathbb C$ is given by
    \begin{equation}\label{eqn:F_pG_ETA_Inverse}
      F_p^{G^{-1}\eta}(z)
      \;:=\;
      \int_{\sqrt{-1}}^z
        (s - s_{-1})^{2}
        (s - s_0)^{-\frac32}
        \prod_{k=1}^{2p}(s - s_k)^{(-1)^{k}\frac12}
      \,ds.
    \end{equation}
We write
\[
  R^p_{G\eta}=R^{\mathrm{ess}}_{X^p_{G\eta}}, 
  \quad
  R^p_{G^{-1}\eta}=R^{\mathrm{ess}}_{X^p_{G^{-1}\eta}},
  \quad
  \omega^p_{G\eta}=\omega_{X^p_{G\eta}},
  \quad
  \omega^p_{G^{-1}\eta}=\omega_{X^p_{G^{-1}\eta}}.
\]
By Remark~\ref{remark: regarding cone angle and order of zero of generalized orthodisk}, the divisors are
\begin{align}
(\omega^p_{G\eta})
  &= P_{t^{\pm}_{-1}}^0\,P_{t_0}^0\,P_{t_1}^2\,P_{t_2}^0\,P_{t_3}^2\,\cdots\,P_{t_{2p}}^0\,P_\infty^{-2}, 
  \label{equation:divisor of gn for p genus}\\
(\omega^p_{G^{-1}\eta})
  &= P_{s^{\pm}_{-1}}^2\,P_{s_0}^{-2}\,P_{s_1}^0\,P_{s_2}^2\,P_{s_3}^0\,\cdots\,P_{s_{2p}}^2\,P_\infty^{-4}.
  \label{equation:divisor of g-1n for p genus}
\end{align}
Here each marked point \(r\) in \(T_{G^j\eta}\) has a lift \(P_r\) (and for \(r=t_{-1}\) or \(s_{-1}\), two lifts \(P_{r^\pm}\)).  One checks that at \(P_{t_0}\) and \(P_\infty\) the divisor of \(\omega^p_{G\eta}\) matches that of \(G_0\,\eta_0\), and at \(P_{s_0}\) and \(P_\infty\) the divisor of \(\omega^p_{G^{-1}\eta}\) matches that of \(G_0^{-1}\,\eta_0\). 

Therefore, a pair of e-generalized orthodisks has been identified such that the resulting surface (if it exists, as in Theorem~\ref{theorem:existence of minimal surface from reflexive orthodisk}) would have an Enneper end and a catenoid end. Such a pair is called the \textit{e-generalized orthodisks of genus $p$ for the Angel surface}. For brevity, this pair will be referred to as the e-generalized orthodisks of genus $p$ in the remainder of the discussion.

Consequently, if an e-reflexive pair of such generalized orthodisks exists, then by Theorem~\ref{theorem:existence of minimal surface from reflexive orthodisk}, such a minimal surface exists  with Weierstrass data
\[
  \Bigl(R^p_{G\eta}\setminus\{P_{t_0},P_\infty\},\,
    G=\tfrac{\omega^p_{G\eta}}{\eta_{X^p_{G\eta}}},\,
    \eta_{X^p_{G\eta}}
  \Bigr).
\]

The remaining task is to search for the pair $(X^p_{G\eta},\,X^p_{G^{-1}\eta})$ that is e-reflexive—that is, conjugate, with a conformal self map of \(\H\cup\R\) takes the marked points to the marked points respecting the order. The idea is that, in the next section, a space of pairs of e-conjugate generalized orthodisks of genus $p$ will be set up, and in the subsequent sections, an e-reflexive pair will be found within that space.

\section{E-conjugate orthodisks with partial
symmetry}
\label{sec:ecs_space}
To  search for a reflexive pair, this section introduces the notion of partial symmetry and defines the pair $(Q_1,\, Q_2)$ for each genus, referred to as a partial symmetric polygon of genus $p \geq 2$. Later, it will be shown that these partial symmetric polygons correspond to the e-conjugate generalized orthodisks of genus $p$. The discussion begins with the definition of two distinct types of ``staircase'' and ``partially symmetric polygons.''
\subsection{Staircases and partially symmetric polygons}\label{subsec:staircases}

We begin by introducing two families of polygonal arcs—\emph{staircases of type~I} and
\emph{staircases of type~II}—which will serve as the building blocks for the
partially symmetric polygons used later in the construction.

\begin{definition}[Type~I elementary blocks]\label{defn:type~I_elementary_block}
	A \emph{Type~I elementary block} is a pair of Euclidean line segments that meet at a common endpoint and are perpendicular to each other making a convex angle as illustrated in the top row of Figure~\ref{fig:elementaryblocks}.  
	The two segments have positive Euclidean lengths \((\ell_1,\ell_2)\).  
	Type~I blocks are classified into three types according to the relation between these lengths. Set of all Type I elementary blocks is denoted as \(\mathcal{I}.\)
\end{definition}

\begin{definition}[Symmetric staircase of type~I of genus~$p$]\label{def:staircase_type_I}
	A \textit{symmetric staircase of type~I} of genus~$p$ (\(p\,\geq\,2\)) is an embedded polygonal arc
	\[
	S=\bigcup_{k=1}^{p-1} B_k\subset\C,
	\]
	where each $B_k$ lies in \(\mathcal{I}.\) The blocks are attached head–to–tail so that the following conditions hold:
	\begin{enumerate}
		\item[\textnormal{(1)}] \textbf{Angle pattern.}  
		As one traverses $S$ starting from the initial point of the leftmost maximal horizontal side, the interior angles (\(i.e\) angles on the left side of the arc) at consecutive vertices alternate between 
		$\pi/2$ (convex) and $3\pi/2$ (concave), beginning with a convex angle at the initial vertex.
		
		\item[\textnormal{(2)}] \textbf{Side lengths.}  
		Side lengths of consecutive blocks may differ arbitrarily;
		
		\item[\textnormal{(3)}] \textbf{Symmetry.}  
		There exists a point $q\in S$ such that $S$ is invariant under the reflection 
		\(\rho:(x,y)\mapsto(-y,-x)\) fixing $q$.
	\end{enumerate}	
\end{definition}
Thus, for a Symmetric staircase of type I of genus \(p,\) \(S=\bigcup\limits_{k=1}^{p-1}B_k,\) the sequence $(B_1,\,\ldots,\,B_{p-1})$ may involve arbitrary choices from $\mathcal{I}$, 
with side lengths chosen independently for each block upto the symmetry condition. See Figure~\ref{fig:staircase_type1_genus3}.
\begin{figure}[ht]
	\centering
	\begin{tikzpicture}[line width=1pt, blue, scale=1.0]
		
		\begin{scope}[shift={(0,0)}]
			\draw (0,0) -- (1.2,0) -- (1.2,2.5);
			\node[below=4pt] at (0.6,0) {\small $I_1$};
		\end{scope}
		
		\begin{scope}[shift={(4,0)}]
			\draw (0,0) -- (3.2,0) -- (3.2,2.5);
			\node[below=4pt] at (1.6,0) {\small $I_2$};
		\end{scope}
		
		\begin{scope}[shift={(9,0)}]
			\draw (0,0) -- (3.2,0) -- (3.2,1.2);
			\node[below=4pt] at (1.6,0) {\small $I_3$};
		\end{scope}
		
		\begin{scope}[shift={(0,-4)}]
			\draw (0,0) -- (0,1.2) -- (3.2,1.2);
			\node[below=4pt] at (1.6,0) {\small $II_1$};
		\end{scope}
		
		\begin{scope}[shift={(4,-4)}]
			\draw (0,0) -- (0,2.8) -- (3.2,2.8);
			\node[below=4pt] at (1.6,0) {\small $II_2$};
		\end{scope}
		
		\begin{scope}[shift={(9,-4)}]
			\draw (0,0) -- (0,2.8) -- (1.4,2.8);
			\node[below=4pt] at (0.7,0) {\small $II_3$};
		\end{scope}
		
		\node at (0, 3.2) {\small \textbf{Type I elementary blocks:}};
		\node at (0, -0.8) {\small \textbf{Type II elementary blocks:}};
		
	\end{tikzpicture}
	
	\caption{Elementary blocks used to construct symmetric staircases.}
	\label{fig:elementaryblocks}
\end{figure}
In a similar spirit to Definition
\ref{defn:type~I_elementary_block}, define \textit{Type~II elementary blocks.}
\begin{definition}[Type~II elementary blocks]
	A \emph{Type~II elementary block} is a pair of Euclidean line segments that meet at a common endpoint and are perpendicular to each other, forming a concave corner, as illustrated in the bottom row of Figure~\ref{fig:elementaryblocks}. The two segments have positive Euclidean lengths \((\ell_1,\ell_2)\).  
	Type~II blocks are classified into three types according to the relation between these lengths.
	The set of all Type~II elementary blocks is denoted by \(\mathcal{II}\).
\end{definition}

\begin{definition}[Symmetric staircase of type~II of genus~$p$]\label{def:staircase_type_II}
	A \textit{symmetric staircase of type~II} of genus~$p$ (\(p\,\ge\,2\)) is an embedded polygonal arc
	\[
	S=\bigcup_{k=1}^{p-1} B_k\subset\C,
	\]
	where each $B_k$ lies in $\mathcal{II}.$ The blocks are attached head–to–tail so that the following conditions hold:
	\begin{enumerate}
		\item[\textnormal{(1)}] \textbf{Angle pattern.}  
		As one traverses $S$ starting from the initial point of the leftmost maximal vertical side, the interior angles at consecutive vertices alternate between 
		$3\pi/2$ (concave) and $\pi/2$ (convex), beginning with a concave angle at the initial vertex.
		
		\item[\textnormal{(2)}] \textbf{Side lengths.}  
		Side lengths of consecutive blocks may differ arbitrarily.
		
		\item[\textnormal{(3)}] \textbf{Symmetry.}  
		There exists a point $q\in S$ such that $S$ is invariant under the reflection 
		\(\rho:(x,y)\mapsto(-y,-x)\) fixing $q$.
	\end{enumerate}
\end{definition}
Thus, for a Symmetric staircase of type II of genus \(p,\) \(S\,=\,\bigcup\limits_{k=1}^{p-1}B_k,\) the sequence $(B_1,\,\ldots,\,B_{p-1})$ may involve arbitrary choices from $\mathcal{II}$, 
with side lengths chosen independently for each block upto the symmetry condition. See Figure~\ref{fig:genus_insertion_P=3}.

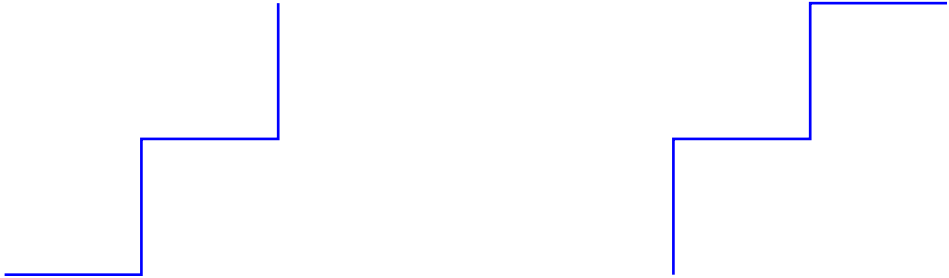
\begin{figure}[ht]
    \centering
    \begin{subfigure}[t]{0.45\textwidth}
        \centering
        \begin{tikzpicture}[scale=0.6, line width=1pt, blue]
            \draw (0,0) -- (2,0) -- (2,3) -- (5,3) -- (5,5);
        \end{tikzpicture}
        \caption{A type I staircase with genus~$p=3$.}
        \label{fig:staircase_type1_genus3}
    \end{subfigure}
    \hfill
    \begin{subfigure}[t]{0.45\textwidth}
        \centering
        \begin{tikzpicture}[scale=0.6, line width=1pt, blue]
            \draw (0,0) -- (0,3) -- (1,3) -- (1,4) -- (4,4);
        \end{tikzpicture}
        \caption{A type II staircase with genus~$p=3$.}
        \label{fig:genus_insertion_P=3}
    \end{subfigure}
    \caption{Comparison of type I and type II staircases of genus~$p=3$.}
    \label{fig:staircase_comparison}
\end{figure}
Next, we introduce two kinds of \emph{partially symmetric polygons}, one of each type, for every genus \(p\ge 2\).

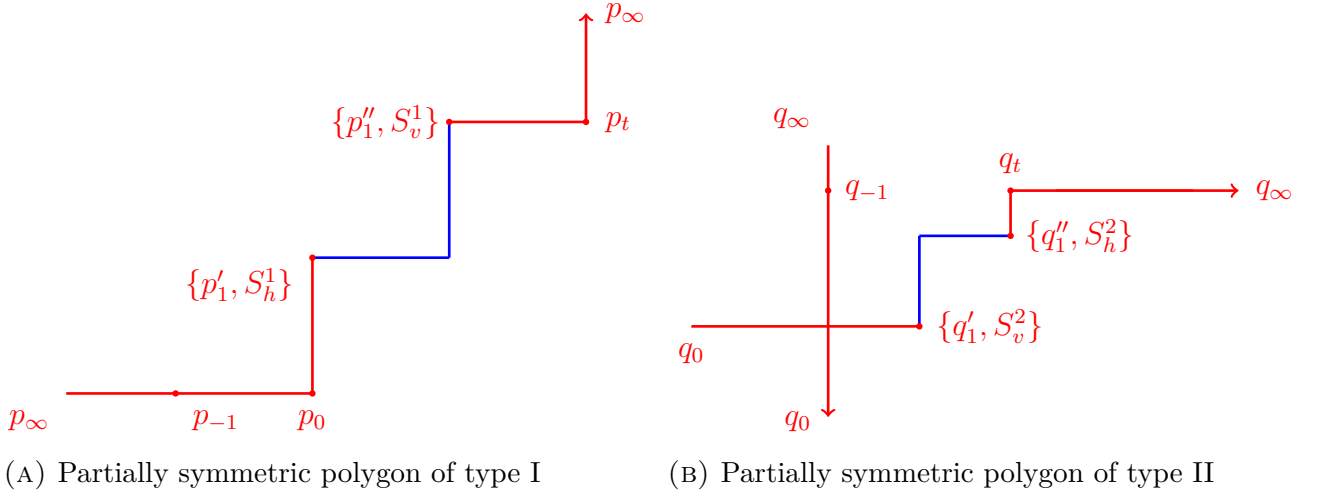
\begin{figure}[ht]
    \centering
    \begin{subfigure}[b]{0.45\textwidth}
        \centering
        \begin{tikzpicture}[line width=1pt, scale=1.2]
            \draw[blue] (1.5,1.5) -- (3,1.5);     
            \draw[blue] (3,1.5) -- (3,3);         

            \draw[red] (0,0) -- (1.5,0);         
            \draw[red] (1.5,0) -- (1.5,1.5);     
            \draw[red] (0,0) -- (-1.2,0);         
            \draw[red] (3,3) -- (4.5,3);
            \draw[red,->] (4.5,3) -- (4.5,4.2);       

            \fill[red] (0,0) circle (1pt) node[below right=2pt] {$p_{-1}$};
            \fill[red] (-1.2,0) node[below left=2pt] {$p_\infty$};
            \fill[red] (1.5,0) circle (1pt) node[below=2pt] {$p_{0}$};
            \fill[red] (1.5,1.5) circle (1pt);
            \fill[red] (4.5,3) circle (1pt) node[right=3pt] {$p_t$};
            \fill[red] (3,3) circle (1pt);
            \fill[red] (4.5,4.2) node[right=3pt] {$p_\infty$};

            \node[red] at (.7,1.2) {$\{p'_1, S^I_h\}$};     
            \node[red] at (2.3,3) {$\{p''_1, S^I_v\}$};    
        \end{tikzpicture}
        \caption{Partially symmetric polygon of type I}
        \label{fig:polygon_type1}
    \end{subfigure}
    \hfill
    \begin{subfigure}[b]{0.45\textwidth}
        \centering
        \begin{tikzpicture}[line width=1pt, scale=1.2]
            \draw[<-,red] (0,-1) -- (0,2);          
            \draw[red] (-1.5,0) -- (1,0);       
            \draw[red] (2.5,1.5) -- (4,1.5);    
            \draw[red] (2,1) -- (2,1.5);   
            \draw[->, red] (2,1.5) -- (4.5,1.5); 

            \draw[blue] (1,0) -- (1,1);         
            \draw[blue] (1,1) -- (2,1);       

            \fill[red] (0,-.75) node[below left=2pt] {$q_0$};
            \fill[red] (0,2) node[above left=2pt] {$q_\infty$};
            \fill[red] (0,1.5) circle (1pt) node[right=2pt] {\(q_{-1}\)};
            \fill[red] (-1.5,0) node[below=2pt] {$q_0$};
            \fill[red] (1,0) circle (1pt) node[right=2pt] {$\{q'_1, S^{II}_v\}$};
            \fill[red] (2,1) circle (1pt) node[right=1pt] {$\{q''_1, S^{II}_h\}$};
            \fill[red] (2,1.5) circle (1pt) node[above=2pt] {$q_t$};
            \fill[red] (4.5,1.5) node[right=2pt] {$q_\infty$};
        \end{tikzpicture}
        \caption{Partially symmetric polygon of type II}
        \label{fig:polygon_type2}
    \end{subfigure}
    \caption{Partially symmetric polygons in the genus 2 case.}
    \label{fig:partially_symmetric_polygons}
\end{figure}

\begin{definition}[Partially symmetric polygon of type I for genus \(p\)]
\label{definition:partially symmetric polygon 1}
Start with the genus-\(1\) generalized orthodisk image as in Figure \ref{fig:genus 1 orthodisk}. 
Take the polygon associated to the differential \(G_{0}\eta_{0}\) and cut it across the vertex \(p_{1}\) along a small line segment passing through \(p_1\) and disjoint from the rest of the polygon in \(\C\).  
This divides the boundary into two piece-wise \(C^{1}\)-arcs
\[
p_{\infty}\,p_{-1}\,p_{0}\,p_{1}^{\prime}
\quad\text{and}\quad
p_{1}^{\prime\prime}\,p_{t}\,p_{\infty},
\]
whose finite endpoints are \(p_{1}^{\prime}\) and \(p_{1}^{\prime\prime}\). 

Choose a genus-\(p\) symmetric staircase of type I, denoted by \(S^I\), with horizontal end \(S^I_{h}\) and vertical end \(S^I_{v}\).  
Identify \(p_{1}^{\prime}\) with \(S^I_{h}\) and \(p_{1}^{\prime\prime}\) with \(S^I_{v}\) so that the interior angles at the glued points 
\(\{p_{1}^{\prime},S^I_{h}\}\) and \(\{p_{1}^{\prime\prime},S^I_{v}\}\) are each \(3\pi/2\).

The region on the left side of the resulting curve in \(\mathbb{C}\) with the curve itself, is called a \textit{partially symmetric polygon of type I}. See Figure \ref{fig:polygon_type1}.
\end{definition}
\begin{definition}[Partially symmetric polygon of type II for genus \(p\)]\label{definition:partially symmetric polygon 2}
Start with the genus-\(1\) generalized orthodisk image as in Figure \ref{fig:genus 1 orthodisk}. 
Take the polygon associated to the differential \(G_{0}^{-1}\eta_{0}\) and cut it across the vertex \(q_{1}\) along a small line segment passing through \(q_1\) and disjoint from the rest of the polygon in \(\C\).  
This divides the boundary into two piece-wise \(C^{1}\)-arcs
\[
q_{\infty}\,q_{-1}\,q_{0}\,q_{1}^{\prime}
\quad\text{and}\quad
q_{1}^{\prime\prime}\,q_{t}\,q_{\infty},
\]
whose finite endpoints are \(q_{1}^{\prime}\) and \(q_{1}^{\prime\prime}\). 

Choose a genus-\(p\) symmetric staircase of type II, denoted by \(S^{II}\), with horizontal end \(S^{II}_{h}\) and vertical end \(S^{II}_{v}\).  
Identify \(q_{1}^{\prime}\) with \(S^{II}_{h}\) and \(q_{1}^{\prime\prime}\) with \(S^{II}_{v}\) so that the interior angles at the glued points 
\(\{q_{1}^{\prime},S^{II}_{h}\}\) and \(\{q_{1}^{\prime\prime},S^{II}_{v}\}\) are each \(\pi/2\).

The region on the left side of the resulting curve in \(\mathbb{C}\) with the curve itself, is called a \textit{partially symmetric polygon of type II}. See Figure \ref{fig:polygon_type2}.
\end{definition}

\begin{definition}[Partially symmetric pair of polygons of genus \(p\)]\label{definition:partially symmetric pair}
Let $Q_1$ be a partially symmetric polygon of type I and $Q_2$ a partially symmetric polygon of type II, both are of genus $p$. Define the genus-$p$ staircases

$$
S^{I} = H^{1}_{1} \cup V^{1}_{1} \cup \cdots \cup H^{1}_{p-1} \cup V^{1}_{p-1},
\qquad
S^{II} = V^{2}_{1} \cup H^{2}_{1} \cup \cdots \cup V^{2}_{p-1} \cup H^{2}_{p-1}
$$
which appear in $Q_1$ and $Q_2$, respectively. Here, for $\varepsilon \in \{1,2\}$, the segment $H^{\varepsilon}_j$ (respectively \ $V^{\varepsilon}_j$) denotes the $j$-th maximal horizontal (respectively \ vertical) edge encountered while traversing the staircase—starting from the horizontal end and moving toward the vertical end in $S^I$, and from the vertical end to the horizontal end in $S^{II}$.

An \textit{ordered pair of partially symmetric polygons of genus $p$} is a pair $(Q_1,\, Q_2)$ such that, after rotating both $Q_1$ and $Q_2$ clockwise by $\pi/4$, the corresponding rotated vectors ${}^r H^\varepsilon_j$ (corresponding to the maximal horizontal edge $H^\varepsilon_j$) and ${}^r V^\varepsilon_j$ (corresponding to the maximal vertical edge $V^\varepsilon_j$) satisfy:
$$
^rH^{1}_{j} = \overline{^rV^{2}_{j}}, 
\qquad
^rV^{1}_{j} = \overline{^rH^{2}_{j}} 
\qquad \text{for every } j = 1, \dots, p-1.
$$
\end{definition}

\subsection{A pair of generalized orthodisks from a pair of partially symmetric polygons}\label{subsec:from polygon to orthodisk} This subsection deals  with two distinct, but ultimately interlinked, objects:

\begin{itemize}
  \item a pair of generalized orthodisks 
        \(\bigl(X_{G\eta}^{p},\,X_{G^{-1}\eta}^{p}\bigr)\)
        defined in Section \ref{section:formalWeierstrassdataForGenusP};
  \item a pair of genus-\(p\) partially symmetric polygons
        \((Q_{1},Q_{2})\) as in
        Definition~\ref{definition:partially symmetric pair}.
\end{itemize}
Starting with a  pair \((Q_{1},\,Q_{2})\), the goal is to construct real parameters so that the corresponding orthodisks
\[
  X_{G\eta}^{p}=X_{G\eta}^{p}(t_{-1},\,t_{0},\,t_{1},\,\dots,\,t_{2p})
  \quad\text{and}\quad
  X_{G^{-1}\eta}^{p}=X_{G^{-1}\eta}^{p}(s_{-1},\,s_{0},\,s_{1},\,\dots,\,s_{2p})
\]
map, via
\eqref{eqn:FpG_ETA} and \eqref{eqn:F_pG_ETA_Inverse}, to the prescribed polygons
\(Q_{1}\) and \(Q_{2}\) (up to rigid motions of~\(\mathbb{C}\)).

To keep the exposition readable, we treat in detail the case \(p=2\).  
The general case is entirely analogous, but involves longer strings of parameters.

\subsubsection{The genus-\(2\) construction}\label{sec:Genus2orthodisk_to_schwarzCristofell}

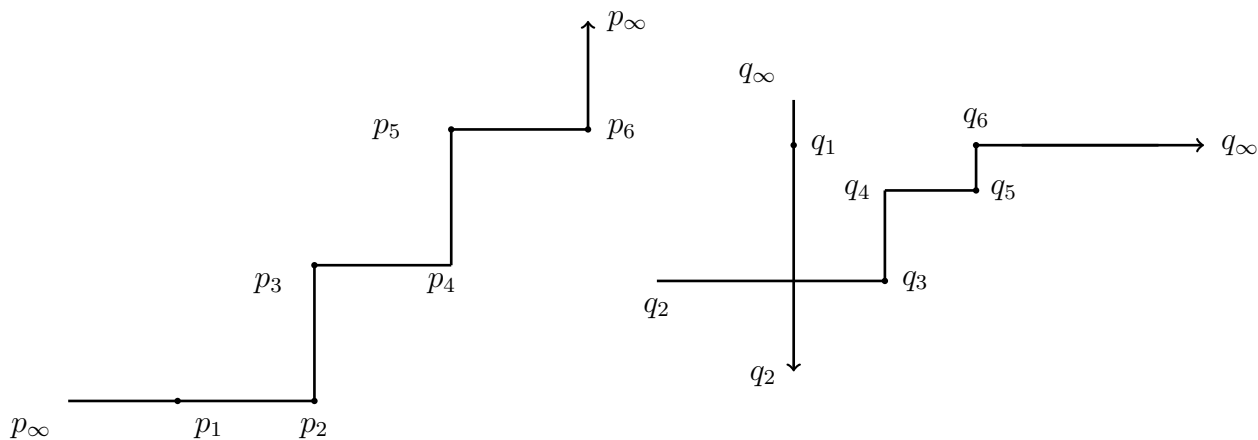
\begin{figure}[ht]
    \centering
    \begin{minipage}{0.48\textwidth}
        \centering
        \begin{tikzpicture}[line width=1pt, scale=1.2]
            \draw[black] (1.5,1.5) -- (3,1.5);     
            \draw[black] (3,1.5) -- (3,3);         

            \draw[black] (0,0) -- (1.5,0);         
            \draw[black] (1.5,0) -- (1.5,1.5);     
            \draw[black] (0,0) -- (-1.2,0);         
            \draw[black] (3,3) -- (4.5,3);
            \draw[black,->] (4.5,3) -- (4.5,4.2);       

            \fill[black] (0,0) circle (1pt) node[below right=2pt] {$p_{1}$};
            \fill[black] (-1.2,0) node[below left=2pt] {$p_\infty$};
            \fill[black] (1.5,0) circle (1pt) node[below=2pt] {$p_{2}$};
            \fill[black] (1.5,1.5) circle (1pt);
            \fill[black] (4.5,3) circle (1pt) node[right=3pt] {$p_{6}$};
            \fill[black] (3,3) circle (1pt);
            \fill[black] (4.5,4.2) node[right=3pt] {$p_\infty$};

            \node[black] at (1,1.3) {$p_{3}$};     
            \node[black] at (2.9,1.3) {$p_{4}$};
            \node[black] at (2.3,3) {$p_{5}$};    
        \end{tikzpicture}
    \end{minipage}
    \hfill
    \begin{minipage}{0.48\textwidth}
        \centering
        \begin{tikzpicture}[line width=1pt, scale=1.2]
            \draw[<-,black] (0,-1) -- (0,2);          
            \draw[black] (-1.5,0) -- (1,0);       
            \draw[black] (2.5,1.5) -- (4,1.5);    
            \draw[black] (2,1) -- (2,1.5);   
            \draw[->, black] (2,1.5) -- (4.5,1.5);

            \draw[black] (1,0) -- (1,1);         
            \draw[black] (1,1) -- (2,1);       

            \fill[black] (0,-.75) node[below left=2pt] {$q_2$};
            \fill[black] (0,2) node[above left=2pt] {$q_\infty$};
            \fill[black] (0,1.5) circle (1pt) node[right=2pt] {$q_{1}$};
            \fill[black] (-1.5,0) node[below=2pt] {$q_2$};
            \fill[black] (1,0) circle (1pt) node[right=2pt] {$q_{3}$};
            \fill[black] (2,1) circle (1pt) node[right=1pt] {$q_{5}$};
            \fill[black] (2,1.5) circle (1pt) node[above=2pt] {$q_{6}$};
            \fill[black] (4.5,1.5) node[right=2pt] {$q_\infty$};
            \node[black] at (0.7,1) {\(q_{4}\)};
        \end{tikzpicture}
            \end{minipage}
    \caption{Image of a pair of generalized orthodisk of genus $2$ : Partially symmetric polygons of type I (left) and type II (right) of genus $2.$}
    \label{fig:genus two orthodisk}
\end{figure}

Let \((Q_{1},Q_{2})\) be a pair of genus-\(2\) partially symmetric polygons.
Denote the finite vertices of \(Q_{1}\) by
\[
  p_{1},\,p_{2},\,p_{3},\,p_{4},\,p_{5},\,p_{6},
\]
with corresponding edge lengths
\(l_{1},l_{2},\dots,l_{5}\).
The interior angles at the finite vertices alternate between
\(\pi/2\) and \(3\pi/2\), starting with \(\pi/2\) at \(p_{2}\). At \(p_{1}\) the angle is~\(\pi\), and at the point at infinity the angle is \(-\pi/2\). See  Figure \ref{fig:genus two orthodisk}.

It should be clarified that $l_1,\,l_2$, and $l_5$  are the same as those of the genus $1$ generalized orthodisk for $G\eta$.
It is a fact (See \cite{stein2010complex}) that a Euclidean polygon determines, up to post-composition by M\"obius transformations, a Schwarz-Christoffel mapping
\[
  F_{0}\,\colon\,\mathbb{H}\longrightarrow Q_{1}
\]
that sends the extended real line \(\mathbb{R}\cup\{\infty\}\) continuously onto the boundary of the polygon.
Because the interior angles -- and hence the exponents in the integrand -- are known, there exist real numbers
\[
  t_{-1}<t_{0}<t_{1}<t_{2}<t_{3}<t_{4}
\]
such that for
\begin{equation}\label{eqn:schwarz_christoffel_F0}
  F_{0}(z)=\int_{\sqrt{-1}}^{z}
    (t-t_{-1})^{0}\,
    (t-t_{0})^{-1/2}\,
    \prod_{k=1}^{4}(t-t_{k})^{(-1)^{k+1}/2}\,dt,
\end{equation}
\(F_0(t_{k})\,=\, p_{k+2}\) for \(k=-1,\,0,\,\ldots,\,4\).

From~\eqref{eqn:schwarz_christoffel_F0} and the prescribed exponents, we obtain a generalized orthodisk as 
\[
  {}^{1}\!T^{2}=\{t_{-1},\,t_{0},\,t_{1},\,t_{2},\,t_{3},\,t_{4}\},
  \qquad
  {}^{1}\!T^{2}_{0}=\{t_0,\,t_{1},\,t_{2},\,t_{3},\,t_{4}\},\,\,
A =  \, \left(1,\, \tfrac12,\, \tfrac32,\,\tfrac12,\,\tfrac32,\,\tfrac12\right).
\]
Similarly, the generalized orthodisk corresponding to the second polygon is obtained.  It is given by 
\[
  {}^{2}\!T^{2}=\{s_{-1},\,s_{0},\,s_{1},\,s_{2},\,s_{3},\,s_{4}\},
  \qquad
  {}^{2}\!T^{2}_{0}=\{s_0,\,s_{1},\,s_{2},\,s_{3},\,s_{4}\},\,\, B=\left(\,3,\,-\tfrac12,\,\tfrac12,\,\tfrac32,\,\tfrac12,\,\tfrac32\right).
\]

The construction used for genus $2$ extends directly to any genus $p \geq 2$ partially symmetric pair $(Q_1, Q_2)$  to get a pair of enhanced generalized orthodisks  of genus $p$ (might not be unique, as one might get more than one conformal polygon with the same image in $\C$). 
\subsection{Partial symmetric polygonal pair to e-conjugate orthodisk}
\label{sec:partialsymmetry_to_conjugacy}
 Let  $(Q_1,\,Q_2)$ be a pair of genus-$p$ partially symmetric polygons as in Definition \ref{definition:partially symmetric pair}. Rotate clockwise both $Q_1$ and $Q_2$ by angle $\pi/4$,  and denote the  new pair by $(Q_1^\prime,\, Q_2^\prime)$.  In view of  Section~\ref{subsec:from polygon to orthodisk},  the pair $(Q_1^\prime,\, Q_2^\prime)$ produces a pair of enhanced generalized orthodisks of genus $p$.
 
 In particular, as described in Section~\ref{subsec:from polygon to orthodisk} (cf. \eqref{eqn:FpG_ETA} and \eqref{eqn:F_pG_ETA_Inverse}), there exist real parameters
$t_{-1} \,<\, t_0 \,<\, t_1 \,<\, t_2 \,<\, \cdots \,<\, t_{2p}$
and
$s_{-1} \,<\, s_0 \,<\, s_1 \,<\, s_2 \,<\, \cdots \,<\, s_{2p},$
such that the corresponding Schwarz-Christoffel map $F_1$, $F_2$ respectively, send these ordered parameters to the vertices of $Q_1^\prime$ and $Q_2^\prime$, respectively. These orthodisks are denoted by $X^p_1$ and $X^p_2$, respectively.

Label the vertices of $Q_1^\prime$ by $P_{-1},\, P_0,\, P_1,\,\dots,\,P_{2p}$, where
$ F_1(t_k)=P_{k}\ (-1\le k\le 2p).$
Likewise, the vertices of $Q_2^\prime$ are labeled $S_{-1},\, S_0,\, S_1,\,\dots,\,S_{2p}$ with
$F_2(s_k)=S_{k}\ (-1\le k\le 2p).$ Note that, what were originally the horizontal and vertical edges have been rotated: the former now lies at an angle of $3\pi/4$
 to the real axis, and the latter at an angle of $\pi/4$. Writing ${}^rH^1_j$ (respectively ${}^rV^1_j$) for the complex displacement of the rotated $j$-th horizontal (respectively vertical) edge  of $Q_1^\prime$, and ${}^rH^2_j,\, {}^rV^2_j$ for those of $Q_2^\prime$, we have:
${}^r H^1_j = \overline{{}^r V^2_j}$ and ${}^r V^1_j = \overline{{}^r H^2_j}$ for  $j=1,\dots,p-1$ (see Definition \ref{definition:partially symmetric pair}).

Let $R^{\mathrm{ess}}_{X^p_1}$, $R^{\mathrm{ess}}_{X^p_2}$ be the hyperelliptic Riemann surfaces obtained from $X^p_1$ and $X^p_2$, and meromorphic 1-forms $\omega_{X^p_i}$, defined as the pulled back of $dz$ under $F_i$ (cf. Remark~\ref{rem:omegaX}). $R^{\mathrm{ess}}_{X^p_1}$, $R^{\mathrm{ess}}_{X^p_2}$ are realized as two-sheeted branched covers of the Riemann sphere with branch points at  $t_0,\, t_1,\,\dots,\,t_{2p}$~ and $ s_0,\, s_1,\,\dots,\,s_{2p}$, respectively.  A canonical homology basis $\{A_j, B_j\}_{j=0}^{p-1}$ is chosen on each $R^{\mathrm{ess}}_{X^p_i}$ as in Subsection \ref{subsection:homology basis for Riemann 
surface} and Remark \ref{rem:homology_basis_notation}.

\begin{figure}[h]
  \centering
  \begin{minipage}{0.45\linewidth}
    \centering
    \def\svgwidth{\linewidth}
    \input{gn.pdf_tex}
    \subcaption{Half arcs corresponding to the homology basis of $R_{X^2_1}^{\mathrm{ess}}$ }
  \end{minipage}
  \hfill
  \begin{minipage}{0.45\linewidth}
    \centering
    \def\svgwidth{\linewidth}
    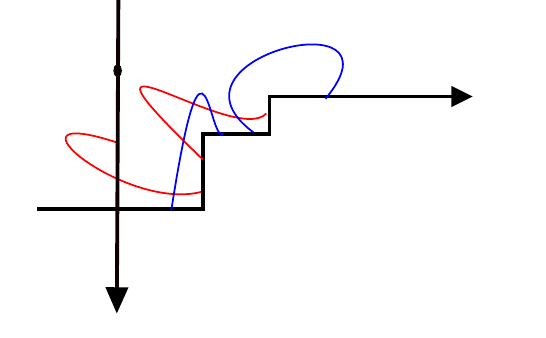
    \subcaption{Half arcs corresponding to the homology basis of $R_{X^2_2}^{\mathrm{ess}}$}
  \end{minipage}
  \caption{Representation of homology basis in the genus \(2\) case.}
    \label{fig:genus2-orthodisk-homology}
\end{figure}

Concretely, $A_0$ is taken as the cycle encircling the branch cut joining $P_0$ to $P_1$ (denote by $P_{0}P_1$) on $Q_1^\prime$ (and likewise encircling $S_0S_1$ on $Q_2^\prime$), and $B_{p-1}$ encircles the branch cut along the edge $P_{2p-1}P_{2p}$ on $Q_1^\prime$ (and $S_{2p-1}S_{2p}$ on $Q_2^\prime$).  The remaining cycles $A_0,\dots,A_{p-1}$ encircle the other rotated vertical staircase edges of $Q_1^\prime$ (resp.  rotated horizontal edges of $Q_2^\prime$), and $B_0,\dots,B_{p-1}$ encircle the other rotated horizontal edges of $Q_2^\prime$ (resp.  rotated vertical staircase edges of $Q_2^\prime$), in corresponding order. See Figure \ref{fig:genus2-orthodisk-homology}.

We now show that the periods of  $\omega_{X^p_1}$ and $\omega_{X^p_2}$ are conjugate appropriately on all these cycles, which will establish that $(X^p_1,\,X^p_2)$  is a conjugate pair of orthodisks.

First, consider the cycle $B_{p-1}$ around the last finite edge pair of the polygons starting from the vertex corresponding to the catenoid end (i.e $t_0$ for $X^p_1$ and $s_0$ for $X^p_2$). By construction, the edge $P_{2p-1}P_{2p}$ in $Q_1$ and the edge $S_{2p-1}S_{2p}$ in $Q_2$ arise from the glued genus  1 reflexive orthodisk pair. In particular, these two edges are complex conjugates of each other.  Computing the period integrals directly from the Schwarz-Christoffel parameterization, the following holds:
$$
\int_{B_{p-1}} \omega_{X^p_1} \;=\; 2\int_{P_{2p-1}}^{P_{2p}}dz \;=\;2\big(P_{2p}-P_{2p-1}\big),
$$
and similarly
$$
\int_{B_{p-1}} \omega_{X^p_2} \;=\; 2\big(S_{2p}-S_{2p-1}\big).
$$
By the construction in Subsection  \ref{subsec:genus 1 orthodisk construction}, $S_{2p}-S_{2p-1}$ is the complex conjugate of $P_{2p}-P_{2p-1}$.  Therefore 
\begin{equation}\label{eq:Bp-1-period}
\int_{B_{p-1}} \omega_{X_1} \;=\; \overline{\int_{B_{p-1}} \omega_{X_2}}.
\end{equation}

Next, consider the cycle $A_0$ encircling the edge $P_0P_1$ on $X^p_1$ and $S_0S_1$ on $X^p_2$. This case is more delicate: in the first orthodisk $X_1$, the edge $P_0P_1$ is a finite length segment, whereas in $X_2$ the corresponding edge $S_0S_1$  extends out to infinity. As a result, a direct evaluation of $\int_{A_0}\omega_{X^p_2}$ via the Schwarz-Christoffel integral is not straightforward. A rigorous comparison of these periods is given below.

Consider the flat surfaces obtained by developing the polygons via the 1-forms: for $i=1,2$ and any genus $p\ge1$, define
\[Y^i_p := \big(\C\setminus\{\text{zeros and poles of }\omega_{X^p_i}\},\,\omega_{X^p_i}\big).\]
The punctured Riemann sphere $Y^i_p$ is the double of $(\mathbb H\cup \mathbb R,\,  \omega_{X^p_i})$ as in Remark \ref{rem:omegaX}.  In particular, $Y^1_p$ contains the straight segment $P_0P_1$ (of finite length) and $Y^2_p$ contains the ray $S_0S_1$ (of infinite length). Now, by construction of the polygons $(Q_1,\,Q_2)$ and hence rotated on $(Q_1^\prime, Q_2^\prime)$,  the local geometry of  $X^p_1$ near the edge $P_0P_1$ is identical to the local geometry of the genus 1 orthodisk near its corresponding edge, and similarly for $X^p_2$ near $S_0S_1$.

Recall $(X_{G_0\eta_0},\,X_{G^{-1}_0\eta_0})$ is  the fixed genus 1 reflexive orthodisk pair used in the construction (see Section~\ref{subsec:genus 1 orthodisk construction} and Figure~\ref{fig:genus 1 orthodisk}).  $p_{-1},\,p_0,\,p_1,\,p_t$ be the the four vertices of the genus 1 polygon $X_{G_0\eta_0}\subset\C$ and $q_{-1},\,q_0,\,q_1,\,q_t$  be the vertices of $X_{G^{-1}_0\eta_0}$ (so $p_0p_1$ and $q_0q_1$ are corresponding edges). The cone angles at $p_0,\,\,p_1$ in $X_{G_0\eta_0}$ agree with those at $P_0,\,P_1$ in $X_1$, and the flat length $|p_1-p_0|$ equals $|P_1-P_0|$ by our choice of parameters.  Therefore we can find a small neighborhood $R^1_1\subset Y^1_1$ containing the segment $[p_0,p_1]$, and a neighborhood $R^1_p\subset Y^1_p$ containing $[P_0,P_1]$, such that $R^1_1$ and $R^1_p$ are isometric as flat surfaces by some map $\varphi$.   Likewise, there are isometric neighborhoods $R^2_1\subset Y^2_1$ around the segment $[q_0,q_1]$ and $R^2_p\subset Y^2_p$ around $[S_0,S_1]$.

Now, in the flat surface  $Y^1_1$ consider a small arc $\gamma$ in the upper half-plane connecting the two sides of the cut $[0,1]$ (that is, running from just to the left of 0 to just to the right of 1 along a path in the upper half-plane).  Doubling this arc $\gamma$ (i.e., reflecting it in the real axis) produces a closed loop in $Y^1_1$ that winds once around the branch cut $(0,1)$ and no other singularity. Without loss of generality, assume that the double of $\gamma$ is contained in $R_1^1$.  In the flat surface $Y^1_p$, the corresponding doubled arc $\varphi(\gamma)$ winds around the cut $P_0P_1$ exactly once (where $\varphi: R^1_1\to R^1_p$ is the local isometry); thus $\varphi(\gamma)$ represents the cycle homotopic to $A_0$ on $R^{\mathrm{ess}}_{X^p_1}$.  Therefore 
$$
\int_{A_0}\omega_{X^p_1} \;=\; \int_{\varphi(\gamma)}\omega_{X^p_1} \;=\; \int_{\gamma} G_0\eta_0\,.
$$
An identical argument applied to $Y^2_1$ and $Y^2_p$ (using the corresponding local isometry $\psi: R^2_1\to R^2_p$) shows that
$$
\int_{A_0}\omega_{X^p_2} \;=\; \int_{\psi(\gamma)}\omega_{X^p_2} \;=\; \int_{\gamma} G^{-1}_0\eta_0\,,
$$ since $\psi(\gamma)$ homotopic to the cycle $A_0$ on $R^{ess}_{X_2^p}$.  But by Proposition~\ref{proposition:genus 1 Angel surface} (the genus 1 conjugacy of $G_0\eta_0$ and $G^{-1}_0\eta_0$), these two genus 1 integrals are conjugate.  It is concluded that
\begin{equation}\label{eq:A0-period}
\int_{A_0}\omega_{X^p_1} \;=\; \overline{\int_{A_0}\omega_{X^p_2}}.
\end{equation}

The next step is to examine the periods on the remaining cycles arising from the staircase segments. Fix $1 \le j \le p - 1$. By the chosen cycles $A_j$ and $B_j$, the following holds:
$$
\int_{A_j}\omega_{X^p_1} \;=\; 2\,{}^r V^1_j,\qquad 
\int_{B_j}\omega_{X^p_1} \;=\; 2\,{}^r H^1_j,
$$
$$
\int_{A_j}\omega_{X^p_2} \;=\; 2\,{}^r H^2_j,\qquad 
\int_{B_j}\omega_{X^p_2} \;=\; 2\,{}^r V^2_j\,.
$$
Now, using the partial symmetry relations ${}^r H^1_j=\overline{{}^r V^2_j}$ and ${}^r V^1_j=\overline{{}^r H^2_j}$ noted above, we immediately find for each $1\le j\le p-1$:
\begin{equation}\label{eq:AjBj-periods}
\int_{A_j}\omega_{X^p_1} \;=\; 2{}^r V^1_j \;=\; 2\overline{{}^r H^2_j} \;=\; \overline{\int_{A_j}\omega_{X^p_2}},\;\;\text{ and }\;\;
\int_{B_j}\omega_{X^p_1} \;=\; 2{}^rH^1_j \;=\; 2\overline{{}^rV^2_j} \;=\; \overline{\int_{B_j}\omega_{X^p_2}}.
\end{equation}
Combining the special cases \eqref{eq:Bp-1-period} and \eqref{eq:A0-period} with the general relation \eqref{eq:AjBj-periods}, it follows that for every cycle $A_j$ or $B_j$ in the canonical homology basis, the period of $\omega_{X^p_1}$ equals the complex conjugate of the period of $\omega_{X^p_2}$. Hence, $(X^p_1,\, X^p_2)$ forms an e-conjugate pair of generalized orthodisks.

\medskip

We conclude this section by formally defining the notion of an e-conjugate pair of generalized orthodisks with partial symmetry.
\begin{definition}\label{defn:essentially_conjugate_pair_partial_symmetry}
	Let 
	\[
	F_p^{G\eta}\colon \mathbb{H}\longrightarrow \C,
	\qquad 
	F_p^{G^{-1}\eta}\colon \mathbb{H}\longrightarrow \C
	\]
	be the Schwarz--Christoffel maps defined in
	\eqref{eqn:FpG_ETA} and \eqref{eqn:F_pG_ETA_Inverse}, respectively.
	A pair of generalized orthodisks
	\(
	\bigl(X_{G\eta},\,X_{G^{-1}\eta}\bigr)
	\)
	is called an \emph{essentially conjugate pair with partial symmetry}
	if after an anti-clockwise rotation
	\(
	R_{\pi/4}\colon \C\longrightarrow\C
	\)
	 the polygonal pair
	\[
	\bigl(
	R_{\pi/4}\bigl(F_p^{G\eta}(\mathbb{H}\cup\R)\bigr),
	\,
	R_{\pi/4}\bigl(F_p^{G^{-1}\eta}(\mathbb{H}\cup\R)\bigr)
	\bigr)
	\]
	coincides with a partially symmetric polygonal pair
	\((Q_1,Q_2)\) introduced in Definition~\ref{definition:partially symmetric pair}.
\end{definition}

Given the partially symmetric polygonal pair $(Q_1,\,Q_2)$ of Definition \ref{definition:partially symmetric pair}, there are in general many pairs of e-conjugate orthodisks $(X_{G\eta},\,X_{G^{-1}\eta})$. Consider a collection of ``equivalence classes'' of such pairs that satisfy Definition \ref{defn:essentially_conjugate_pair_partial_symmetry}; our goal is to locate within this parameter space an \emph{e-reflexive}  pair.

\section{Space of e-conjugate pair of orthodisks with partial symmetry}
\label{sec:spaceOFessentiallyCOnjugateOrthodiskWPS}

This section defines the moduli space of genus-$p$ orthodisk pairs that are e-conjugate and satisfy the partial symmetry constraints introduced earlier. The construction proceeds in two stages: first, by describing the moduli of the underlying partially symmetric polygonal pairs, and then by incorporating the period (conjugacy) conditions to define the full moduli space of orthodisk pairs.

\subsection{Partially symmetric polygonal pairs and the moduli $\Sigma_p$}\label{subsec:moduli_polygonal_pair}

Recall from Section~\ref{sec:ecs_space} that a partially symmetric polygonal pair of genus $p$ consists of two polygonal curves $(Q_1,\,Q_2)$ in the plane, each built by inserting a symmetric staircase consisting of $p-1$ elementary blocks into a fixed genus 1 polygon in a certain way.  Let $\widetilde{\Sigma}_p$ denote the set of all such partially symmetric polygonal pairs of genus $p$. The (orientation-preserving) Euclidean isometry group $\mathrm{Iso}^+(\mathbb{R}^2)$ acts naturally on $\widetilde{\Sigma}_p$ by diagonally applying rigid motions to the pair $(Q_1,\,Q_2)$. We define the moduli space of partially symmetric polygon pairs of genus $p$, up to rigid motion, as the quotient:

$$
\Sigma_p \;:=\; \widetilde{\Sigma}_p \Big/ \mathrm{Iso}^+(\mathbb{R}^2)\,. 
$$

Intuitively, $\Sigma_p$ records the intrinsic shape parameters of the polygonal pair, forgetting extrinsic position or rotation in the plane. By construction, each element of $\Sigma_p$ encodes the lengths of edges and the angles of the partially symmetric pair $(Q_1,\,Q_2)$, subject to the symmetry constraints.

To describe $\Sigma_p$ more concretely, it is convenient to begin with a base
configuration arising in genus~$1$. Fix a particular genus 1 e-reflexive orthodisk pair $(\tilde X_{G_0\eta_0},\,\tilde X_{G_0^{-1}\eta_0})$ as in Lemma \ref{lemma:existence_of_genus_one_orthodisk}. Let $(Q_1^{(1)},\, Q_2^{(1)})$ be the corresponding polygonal pair of genus $1$. This base configuration is denoted by $\lambda_0 \in \Sigma_1$. All the cone angles in this genus 1 pair are fixed by construction; thus $\lambda_0$ can be uniquely characterized by the lengths of three consecutive edges in one of its polygonal representations.  More concretely, if $p_{-1},\, p_0,\, p_1,\, p_t$ are four consecutive vertices along $Q_1^{(1)}$ (with $p_{-1},\,p_0,p_1$ and $p_t$ finite and $p_{\infty}$ corresponding to a point at infinity), then we set
$$
\ell_{-1} := |p_{-1}-p_0|,\qquad 
\ell_{0} := |p_{0}-p_{1}|,\qquad 
\ell_{t} := |p_{1}-p_{t}|\,,
$$
the three positive lengths. It is possible to identify
$$
\lambda_0 \;\equiv\; (\ell_{-1},\,\ell_{0},\,\ell_{t}) \;\in\; (0,\infty)^3\,,
$$
using these three edge lengths as coordinates for the genus 1 moduli $\Sigma_1$.  In particular, $\lambda_0$ determines both $Q_1^{(1)}$ and $Q_2^{(1)}$, since the genus 1 pair is e-reflexive. 

Now consider an arbitrary genus $p$ partially symmetric pair $(Q_1,\, Q_2) \in \widetilde{\Sigma}_p$. The pair $(Q_1,\, Q_2)$ is said to be obtained from the base $\lambda_0$ by inserting a staircase of length $p - 1$ if $Q_1$ and $Q_2$ reduce to the base polygons $Q_1^{(1)}$, $Q_2^{(1)}$ when the staircase segments are removed. Let $\Delta_{\lambda_0,p} \subset \Sigma_p$ denote the set of all genus $p$ partially symmetric polygon pairs whose underlying genus-1 shape is $\lambda_0$ in this sense. The set $\Delta_{\lambda_0,p}$ may be regarded as the slice of the moduli space $\Sigma_p$ obtained by ``adding $p-1$ handles'' (via staircase inserts) to the fixed base configuration $\lambda_0$.

An element $\lambda\in \Delta_{\lambda_0,p}$ is determined by the additional lengths introduced by the staircase inserts. In particular, each staircase of length $p-1$ introduces $p-1$ new elementary blocks (Figure \ref{fig:elementaryblocks}) on each of $Q_1$ and $Q_2$. Due to the partial symmetry of $Q_1$ and $Q_2$, the lengths of the last $p-1$ edges of staircases of $Q_1$ must coincide with the lengths of the first $p-1$ edges (in reverse order) of $Q_1$. The same holds for $Q_2$. Thus, there are exactly $p-1$ independent edge-length parameters for the genus-$p$ pair beyond those already in $\lambda_0$. We can label these independent lengths by $\ell_1,\,\ell_2,\,\dots,\,\ell_{p-1}$, corresponding to the successive finite edges along (say) $Q_1$ from the first inserted staircase block to the $(p-1)$-th block. Formally, enumerate the new vertices along $Q_1$ as $P_{t_1},\, P_{t_2},\, \dots,\, P_{t_{2p-1}}$ (in order along the boundary of $Q_1$), where $P_{t_1}$ is the first new vertex after the splitting point and $P_{t_{2p-1}}$ is the last new vertex. Then define the vector of edge lengths

$$
J(\lambda) \;:=\; \Big( |P_{t_1}-P_{t_2}|,\; |P_{t_2}-P_{t_3}|,\; \dots,\; |P_{t_{p-1}}-P_{t_p}| \Big) \,\in\, (0,\infty)^{\,p-1}\,. 
$$

By the reflection symmetry of $Q_1$, the remaining $p-1$ edges (from $P_{t_p}$ to $P_{t_{2p-1}}$) have the same lengths, so no further parameters appear. In this way, each $\lambda\in\Delta_{\lambda_0,p}$ is encoded by a $(p-1)$-tuple $J(\lambda) = (\ell_1,\dots,\ell_{p-1})$ of positive real numbers.  This provides a coordinate chart on $\Delta_{\lambda_0,p}$.  This naturally induces a topology on \(\Delta_{\lambda_0,p}\), and it is a homeomorphism onto its image.  We thus identify
\begin{equation}\label{eq:Delta-chart}
   \Delta_{\lambda_0,p}\;\cong\;(0,\infty)^{p-1},
   \qquad \lambda\longleftrightarrow (\ell_1,\,\dots,\,\ell_{p-1}),
\end{equation}
where \(\ell_k = |P_{t_k}-P_{t_{k+1}}|\).

\subsection{Space of e-conjugate orthodisk pairs and the space $T_{\lambda_0,p}$}
\label{subsection:geometry_moduli_space}
Define: 
 
\[
\widetilde{\Lambda}_{\lambda_0,p} := 
\left\lbrace\, 
\begin{aligned}
&\bigl(X^p_{G\eta}(t_{-1}, t_0,\ \ldots\ ,t_{2p}),\; X^p_{G^{-1}\eta}(s_{-1}, s_0,\ \ldots \ ,s_{2p})\bigr)
\;\Big|\; \\
&\quad t_{-1} < t_0 < t_1 < \cdots < t_{2p}, \\
&\quad s_{-1} < s_0 < s_1 < \cdots < s_{2p}, \\
&\quad \text{these define an e-conjugate, partially symmetric} \\
&\quad \text{orthodisk pair of genus $p$, as in Definition~\ref{defn:essentially_conjugate_pair_partial_symmetry}}, \text{ obtained from } \lambda_0
\end{aligned}
\,\right\rbrace.
\]

$\mathrm{Iso}^+(\R^2)$ acts on $\widetilde{\Lambda}_{\lambda_0,p}$ diagonally as action on the corresponding $(Q_1,\,Q_2)$.  Define \[
   T_{\lambda_0,p} 
   \;:=\;
   \widetilde{\Lambda}_{\lambda_0,p}
   \;\big/\;\mathrm{Iso}^+(\R^2).
\]

Each orthodisk pair determines an ordered pair of conformal structures on the $M= \mathbb{CP}^1 \setminus \{r_{-1},r_0,r_1,\dots,r_{2p},\,\infty\}$. Using exactly the same notation as in \cite[Section 4.2]{weber1998minimal}, take \(T_{\lambda_0,p}\subset \mathrm{Teich}(M)\times\mathrm{Teich}(M)\) whose points are equivalence classes of pair of $2p+3$ punctured spheres that arise from a partially symmetric conjugate generalized orthodisk pair of genus \(p\).  
By \cite[Lemma 3.1.4]{weber2002teichmuller}, the map 
\[
   \Pi: T_{\lambda_0,p} \longrightarrow \Delta_{\lambda_0,p},
   \qquad
   \Pi([X^p_{G\eta}, \,X^p_{G^{-1}\eta}]) = (\ell_1,\dots,\ell_{p-1}),
\]
which sends a pair of marked conformal structures to its \(p-1\) independent staircase edge lengths, is a local diffeomorphism.  As in the case of zigzags, these lengths provide a local chart in the moduli space of the e-conjugate pair with partial symmetry.

Note that since for $\lambda_0$, we fixed $(\ell_{-1},\, \ell_0,\, \ell_t)$ and for higher genus $(p\,\geq\, 2)$, we are looking for the variation of the length of steps of stairs only which is similar to changing the length of steps of zigzags as defined in \cite{weber1998minimal}.  Therefore, most of the Teichm\"{u}ller theory aspects of the discussion here turn out to be the same as the case of zigzags discussed by Weber and Wolf in \cite{weber1998minimal}.   

\section{Height function on the space \texorpdfstring{$T_{\lambda_0,p}$}{T_lambda0,p}}\label{subsection:Teich_extremal_length}
This section introduces the height function on the space $T_{\lambda_0,p}$. The height function serves as a real-analytic measure of the failure of e-reflexivity for a given pair, and its vanishing actually characterizes reflexive pairs, which correspond to solutions of the period problem for the minimal surface construction.

Let $(\zeta_1,\, \zeta_2)= [X^p_{G\eta},\, X^p_{G^{-1}\eta}]\,\in\, T_{\lambda_0,p}$.  Fix a standard collection of homotopy classes of simple closed curves  on the underlying punctured sphere $M = \mathbb{CP}^1 \setminus \{r_{-1}, r_0, r_1, \ldots, r_{2p}, \infty\}$, corresponding to  $X^p_{G\eta}$ and $X^p_{G^{-1}\eta}$, denoted by
\begin{equation}\label{eqn:Homotopyclassofloop}
\Gamma_{-1},\, \Gamma_0,\, \Gamma_1,\, \ldots,\, \Gamma_{2p-1},
\end{equation}
where each $\Gamma_j$ encircles precisely the pair of punctures $r_j$ and $r_{j+1}$.

\subsection{Extremal Length \cite[Section 2.3]{weber1998minimal}} 
For each $j \in \{-1,\,0,\,1,\,\dots,2p-1\}$, let $\mathrm{Ext}_{G\eta}(\Gamma_j)$ denote the extremal length of the loops $\Gamma_j$  when calculated in flat metric on punctured sphere corresponding to  $X^p_{G\eta}$ (induced by the meromorphic 1-form $\omega_{X^p_{G\eta}}$ associated to the orthodisk $X^p_{G\eta}$). Likewise, let $\mathrm{Ext}_{G^{-1}\eta}(\Gamma_j)$ be the extremal length of the $\Gamma_j$ on the  flat surface  corresponding to conjugate orthodisk $X^p_{G^{-1}\eta}$. For a given pair $\zeta = (\zeta_1,\, \zeta_2) \in T_{\lambda_0,p}$, this construction yields two collections of positive real numbers:
$$
\big\{\mathrm{Ext}_{G\eta}(\Gamma_j)\big\}_{j=-1}^{2p-1}, \qquad 
\big\{\mathrm{Ext}_{G^{-1}\eta}(\Gamma_j)\big\}_{j=-1}^{2p-1}. 
$$

Henceforth, we adopt the notation \( \mathrm{Ext}_{G\eta}(\Gamma_j; \zeta_1) \) and \( \mathrm{Ext}_{G^{-1}\eta}(\Gamma_j; \zeta_2) \) to emphasize that \( G^{\pm 1}\eta \) varies with the parameters \( \zeta_1 \) and \( \zeta_2 \), respectively. This highlights the dependence of extremal lengths on coordinates \( \zeta_1 \) and \( \zeta_2 \).

\subsection{Height function}\label{sec:height_func}
For any free homotopy class $[\gamma]$ of loops on $M$, define the nonnegative quantity
$$
H^p_{[\gamma]}(\zeta) \;=\; \Big(e^{\,\mathrm{Ext}_{G\eta}([\gamma];\,\zeta)} - e^{\,\mathrm{Ext}_{G^{-1}\eta}([\gamma];\,\zeta)}\Big)^2 \;+\; \Big(e^{\,1/\mathrm{Ext}_{G\eta}([\gamma];\,\zeta)} - e^{\,1/\mathrm{Ext}_{G^{-1}\eta}([\gamma];\,\zeta)}\Big)^2. 
$$

$H^p_{[\gamma]}(\zeta)$ vanishes if and only if the extremal length of the loop $[\gamma]$ is the same in both flat metrics determined by $\zeta$. In particular, consider the specific collection of homotopy classes of loops $\Gamma_j$ encircling adjacent punctures as above (see \eqref{eqn:Homotopyclassofloop}). Define the height function $H_p: T_{\lambda_0,p} \longrightarrow \mathbb{R}_{\ge 0}$ by summing the contributions of these loops:
\begin{equation}\label{eq:height_function}
H_p(\zeta) \;:=\; H^p_{\Gamma_{-1}}(\zeta) \;+\; H^p_{\Gamma_0}(\zeta) \;+\; H^p_{\Gamma_{2p-1}}(\zeta) \;+\; \sum_{j=1}^{p-1} H^p_{\Gamma_j}(\zeta). 
\end{equation}
The function $H_p$ defined above is the same (adapted in the new setup) height function introduced in \cite{weber1998minimal}, and hence $H_p$ is a proper function on $T_{\lambda_0,p}$ \cite[Theorem 4.6.1]{weber1998minimal}.

If the pair of generalized orthodisks is e-reflexive, then it is direct to see that $H^p(\zeta)=0$.  Conversely, if $H^p(\zeta_0)=0,$ it follows that 
    \begin{equation*}
    \mathrm{Ext}_{G\eta}(\Gamma_j;\zeta_0)=\mathrm{Ext}_{G^{-1}\eta}(\Gamma_j;\zeta_0)\text{ for }j=-1,\,0,\,1,\,\dots,\,p-1,\text{ and } 2p-1.\end{equation*}

Due to the symmetry, 
\(\mathrm{Ext}_{G\eta}(\Gamma_j;\zeta_0)=\mathrm{Ext}_{G^{-1}\eta}(\Gamma_j;\zeta_0)\text{ for }j=-1,\,0,\,1,\,\ldots,\,2p-1.\)
In other words, the pair of conformal structure induced by \(\zeta_0\) on \(\C\setminus\{r_{-1},\,r_{0},\,r_1,\ldots,\,r_{2p}\}\) are conformal to each other. Hence, \(\zeta_0\) is reflexive.

Therefore, it suffices to find at least one zero of $H_p$ in $T_{\lambda_0,p}$;  this will be carried out in Section \ref{section:Existence of Essentially Reflexive Generalized Orthodisks of Genus p}.

\section{Tangent space $T_{\texorpdfstring{\zeta}{zeta}}(T_{\texorpdfstring{\lambda_0,p}{lambda0,p}})$}
\label{section:tangent_space}

 This section identifies the tangent space of $T_{\lambda_0,p}$ by adapting the strategy used by Weber-Wolf for the zigzag family in \cite{weber1998minimal}. For every genus $p \ge 2$, the deformation is confined to the staircase portion of the surface, so the analytical framework coincides with the zigzag case. The only procedural difference from the Weber-Wolf setting is that a genus-one orthodisk is first fixed, and the staircase building blocks are subsequently inserted. All discussions therefore carry over verbatim once expressed in this revised notation, and the relevant modifications are recorded below.

As in \cite{weber1998minimal},  consider a class of diffeomorphisms of \(\mathbb{C}\) that map a pair of orthodisks to another such pair, and compute their infinitesimal Beltrami differentials. These differentials will constitute the tangent space at \(\lambda\).

Let the image of $\zeta \in T_{\lambda_0,p}$ in $\mathbb{C} \times \mathbb{C}$ be denoted by $\lambda=(\Omega_{G\eta},\, \Omega_{G^{-1}\eta})\in\Delta_{\lambda_0,p}$ via the map \eqref{eqn:FpG_ETA} and \eqref{eqn:F_pG_ETA_Inverse} respectively. Select an edge $E$ of the type I staircase in the partially symmetric polygon $\Omega_{G\eta}$, and its corresponding edge in the type II polygon $\Omega_{G^{-1}\eta}$, which is also denoted by $E$ by slight abuse of notation. Suppose $E$ is a horizontal edge of $\Omega_{G\eta}$. A diffeomorphism $f^E_\epsilon$ of $\mathbb{C}$ is applied, which infinitesimally displaces the edge $E$, while remaining the identity outside a compact neighborhood of $E$. This deformation breaks the partial symmetry of the staircase in $\Omega_{G\eta}$ along the diagonal line $y = -x$.

\begin{figure}[ht]
    \centering
    \includegraphics[width=0.8\linewidth]{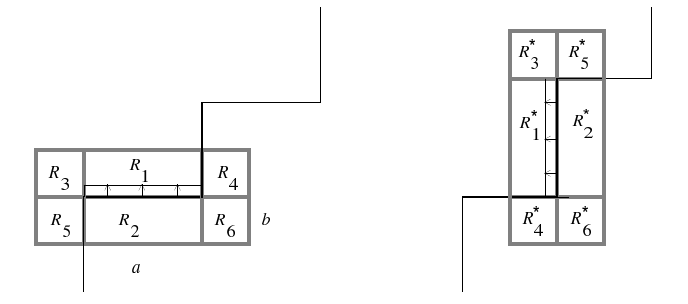}
    \caption{Representation of \(F_{\epsilon}\); copied from \cite[Page 774]{weber2002teichmuller}.}
    \label{fig: pictorial representation of the map taking zigzag to zigzag}
\end{figure}

Let $E^*$ denote the edge that is the reflected image of $E$ in the stair.  To restore this symmetry, we have to adjust $E^*$ accordingly. Introduce a second diffeomorphism \({f^{E^*}_{\epsilon}}\), which infinitesimally displaces \(E^*\). The explicit forms of the maps \(f^E_\epsilon\) and \(f^{E^*}_{\epsilon}\) are exactly those described in equations (5.1a) and (5.1b) of~\cite{weber2002teichmuller}.  In Figure \ref{fig: pictorial representation of the map taking zigzag to zigzag}, these maps are explained as push and pull maps. 

Let \(F_{\epsilon} := f^{E^*}_{\epsilon} \circ f^E_\epsilon\). Then \(F_{\epsilon}\) is a diffeomorphism of \(\mathbb{C}\) that transforms \(\Omega_{G\eta}\) into another partially symmetric object.  We call the map $ F_{\epsilon}$ the ``pushing out and pulling in''  map for the edge $E$.

Let $\nu_\epsilon \,:=\, \frac{(f_\epsilon^E)_{\overline{z}}}{(f_\epsilon^E)_{z}}$ represents the Beltrami differential of $f_\epsilon^E$, and 
define $\dot{\nu} \,= \,\left. \frac{d}{d\epsilon}\right|_{\epsilon=0}\nu_\epsilon$. Similarly, let $\dot{\nu^{*}}$ denote the infinitesimal 
Beltrami differential of $f^{E^*}_{\epsilon}$. Expressions for $\dot{\nu}$ and $\dot{\nu^{*}}$ are given in
\cite[(5.2)(a)]{weber2002teichmuller} and \cite[(5.2)(b)]{weber2002teichmuller} respectively.

We take $\dot{\mu} \,= \,\dot{\nu} + \dot{\nu}^*$. This is a Beltrami differential supported on a bounded domain in
$\mathbb{C}$. This pair of Beltrami differentials lifts to a pair 
\begin{equation}\label{tangentvector}
    \dot{\mu}\,=\,(\dot{\mu}_{\Omega_{G\eta}},\, \dot{\mu}_{\Omega_{G^{-1}\eta}}).
\end{equation}

The above defined $\dot{\mu}$ represents a tangent vector to $T_{\lambda_0,p}$ at $\zeta$ under certain equivalence relation.  The above process will yield different tangent vectors for different ``pushing out and pulling in'' maps.

\section{Existence of e-reflexive generalized orthodisks of genus $p$}\label{section:Existence of Essentially Reflexive Generalized Orthodisks of Genus p}

We now show that for each integer $p\,\ge\, 1$ there exists at least one e-reflexive pair of generalized orthodisks of genus $p$ (in the sense of Definition ~\ref{definition: essentially reflexive orthodisk}).   

The case  of $p=1$ is true by Lemma \ref{lemma:existence_of_genus_one_orthodisk}. To prove the general case by induction, assume:
\begin{assumption}\label{assumption:reflexive_pair_genus_p_minus_1}
    There exists an e-reflexive pair of generalized orthodisks of genus~\( p-1 \). Let us call it $\zeta^{p-1}.$
\end{assumption}

\begin{table}[H]
    \centering
    \renewcommand{\arraystretch}{1.5}
    \begin{tabular}{|p{6cm}|p{6cm}|}
        \hline
        \textbf{Edge pair of the partially symmetric polygon of type I} & \textbf{Edge pair of the partially symmetric polygon of type II} \\
        \hline
        \( (P_{t_1}P_{t_2},\,P_{t_{2p-2}}P_{t_{2p-1}}) \) & \( (P_{s_1}P_{s_2},\,P_{s_{2p-2}}P_{s_{2p-1}}) \) \\
        \hline
        \( (P_{t_2}P_{t_3},\,P_{t_{2p-3}}P_{t_{2p-2}}) \) & \( (P_{s_2}P_{s_3},\,P_{s_{2p-3}}P_{s_{2p-2}}) \) \\
        \hline
        \( \vdots \) & \( \vdots \) \\
        \hline
        \( (P_{t_{p-1}}P_{t_p},\,P_{t_{p}}P_{t_{p+1}}) \) & \( (P_{s_{p-1}}P_{s_p},\,P_{s_p}P_{s_{p+1}}) \) \\
        \hline
    \end{tabular}
    \caption{Correspondence between edge pairs of polygons of types I and II.}
    \label{tab:my_label}
\end{table}
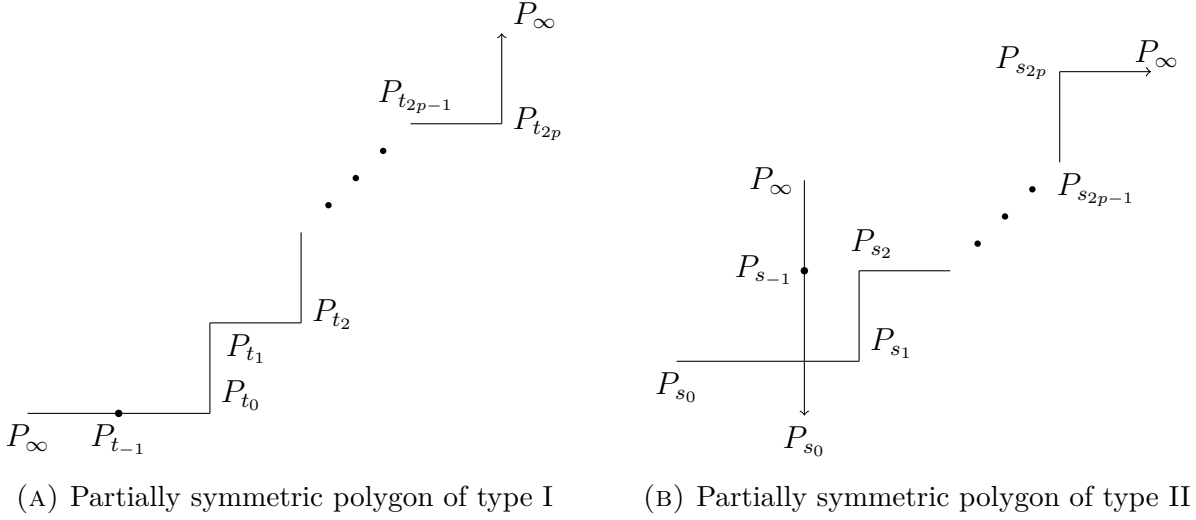
\begin{figure}
\centering
\begin{subfigure}{0.48\textwidth}
\centering
\begin{tikzpicture}[scale=1.2]

\draw (-1,0) -- (0,0);
\node[below] at (-1,0) {$P_{\infty}$};
\draw (0,0) -- (1,0);
\node[below] at (0,0) {$P_{t_{-1}}$};
\draw[fill=black] (0,0) circle (1pt);

\draw (1,0) -- (1,1);
\node[right] at (1,0.2) {$P_{t_0}$};

\draw (1,1) -- (2,1);
\node[below] at (1.4,1) {$P_{t_1}$};

\draw (2,1) -- (2,2);
\node[right] at (2,1.1) {$P_{t_2}$};



\foreach \i in {1,...,3}
  \fill (2 + 0.3*\i, 2 + 0.3*\i) circle (1pt);

\draw (3.2,3.2) -- (4.2,3.2);
\node[left] at (3.8,3.5) {$P_{t_{2p-1}}$};

\draw[->] (4.2,3.2) -- (4.2,4.2);
\node[right] at (4.2,3.2) {$P_{t_{2p}}$};
\node[right] at (4.2,4.4) {$P_{\infty}$};
\end{tikzpicture}
\caption{Partially symmetric polygon of type I}
\end{subfigure}
\hfill
\begin{subfigure}{0.48\textwidth}
\centering
\begin{tikzpicture}[scale=1.2]

\draw (-1,0) -- (1,0);
\node[below] at (-1,0) {$P_{s_0}$};
\draw[->] (0.4,2) -- (0.4,-0.6);
\node[left] at (0.4,2) {$P_{\infty}$};
\node[left] at (0.4,1) {$P_{s_{-1}}$};
\node[below] at (0.4,-0.6) {$P_{s_0}$};
\draw[fill=black] (0.4,1) circle (1pt);

\draw (1,0) -- (1,1);
\node[right] at (1,0.2) {$P_{s_1}$};

\draw (1,1) -- (2,1);
\node[left] at (1.5,1.3) {$P_{s_2}$};




\foreach \i in {1,...,3}
  \fill (2 + 0.3*\i, 1 + 0.3*\i) circle (1pt);

\draw (3.2,2.2) -- (3.2,3.2);
\node[below] at (3.6,2.2) {$P_{s_{2p-1}}$};

\draw[->] (3.2,3.2) --  (4.2,3.2);
\node[left] at (3.2,3.3) {$P_{s_{2p}}$};
\node[right] at (3.9,3.4) {$P_{\infty}$};
\end{tikzpicture}
\caption{Partially symmetric polygon of type II}
\end{subfigure}
\caption{A partially symmetric pair of polygons of genus $p$}
\label{figure:p_genus_partially_symmetric_polygon_pair}
\end{figure}
Recall from Section \ref{sec:spaceOFessentiallyCOnjugateOrthodiskWPS} that the map
\[
\Pi \,\colon\, T_{\lambda_0,p} \;\longrightarrow\; \Delta_{\lambda_0,p}
\]
is a local diffeomorphism. Hence, for each $\zeta \in T_{\lambda_0,p}$, there is an isomorphism of real tangent spaces
\[
d\Pi_\zeta \;:\; 
T_\zeta T_{\lambda_0,p}
\;\xrightarrow{\;\simeq\;}
T_{\Pi(\zeta)= \lambda}\Delta_{\lambda_0,p}
. \]
The real dimension of  $T_\zeta T_{\lambda_0,p}$ is $p-1$ and tangent vectors at \(\zeta\) are realized by a family of Beltrami differentials as discussed in Section \ref{section:tangent_space}.  These are obtained by perturbing the lengths of the edge-pairs (see Figure~\ref{figure:p_genus_partially_symmetric_polygon_pair}) of the type I and type II partially symmetric polygons corresponding to $\lambda$, respectively,  as described in Table~\ref{tab:my_label} .  Perturbing these edges produces a set of infinitesimal Beltrami pairs
\[
(\dot{\mu}^j,\;\dot{\widetilde\mu}^j)
\]
with compact support in small neighborhoods of those edge‐pairs where $j=1,\,2,\,\dots,\,p-1$  (cf. Section \ref{section:tangent_space}). Their equivalence class is denoted by
\[
(\dot{\mu_0}^j,\;\dot{\widetilde\mu}^j_0)
\;\in\;
T_\zeta T_{\lambda_0,p},
\]
where two pairs \((\dot\mu,\dot{\widetilde\mu})\) and \((\dot\mu',\dot{\widetilde\mu}')\) represent the same tangent vector if
\[
\int_{M^{G\eta}_{0,2p+2}} \!\!\Phi_1\,(\dot\mu-\dot\mu')
\;=\;
0
\quad\text{and}\quad
\int_{M^{G^{-1}\eta}_{0,2p+2}} \!\!\Phi_1'\,(\dot{\widetilde\mu}-\dot{\widetilde\mu}')
\;=\;
0,
\]
for every 
\(\Phi_1\in \mathrm{QD}\bigl(M^{G\eta}_{0,2p+2}\bigr)\)
and
\(
\Phi_1' \in \mathrm{QD}\bigl(M^{G^{-1}\eta}_{0,\,2p+2}\bigr).
\)
Here \( M^{G^{j}\eta}_{0,\,2p+2} \) for \( j = \pm 1 \) denotes a pair of Riemann surfaces that are topologically \( \mathbb{C} \) with $2p+2$ punctures, and the conformal structures are induced by \( \zeta \). Also \( \mathrm{QD}(M) \) denotes the space of holomorphic quadratic differentials on the Riemann surface \( M \) \cite[Section 2.2]{weber2002teichmuller}.

\subsection{Variation of extremal length and admissible foliation}
\label{subsec:extremal_length}

Gardiner’s formula gives the first‐order variation of extremal length \cite[Equation 2.2]{weber2002teichmuller}. For any free homotopy class $[\gamma]$, the following holds:
\[
\left(d\operatorname{Ext}_{G\eta}\!\bigl([\gamma];\cdot\right)\big\vert_{\zeta}\,\mu_0,\;
  d\operatorname{Ext}_{G^{-1}\eta}\!\bigl([\gamma];\cdot\bigr)\big\vert_{\zeta}\,\widetilde{\mu}_0
\bigr)
\;=\;
4\bigl(\operatorname{Re}\!\int_{M^{G\eta}_{0,\,2p+2}}\!\Phi_{[\gamma]}\,\mu_0,\;
  \operatorname{Re}\!\int_{M^{G^{-1}\eta}_{0,\,2p+2}}\!\Phi_{[\gamma]}'\,\widetilde{\mu}_0
\bigr).
\]
Here \(\Phi_{[\gamma]}\) is the holomorphic quadratic differential whose horizontal foliation consists of curves that connect the same edges as $\gamma$. See \cite[Subsection 5.3.2]{weber2002teichmuller}. The relevant homotopy classes are $\{\Gamma_j\}_{j=-1}^{2p-1}$, defined in  \eqref{eqn:Homotopyclassofloop}.

Next, we discuss the definition of an admissible foliation \(\Phi_{\Gamma_j}\) for an edge $E$. Recall that \(\Gamma_j\) is the free homotopy class of a topological loop encircling the \(j^{\text{th}}\) and \((j+1)^{\text{th}}\) punctures on the topological sphere with \(2p+3\) punctures. 
\begin{definition}[Admissible foliation for the edge]\cite{weber2002teichmuller}]\label{defn:admissible_edge}
The quadratic differential \(\Phi_{\Gamma_j}\) is admissible for the  edge \(E\) connecting vertices \(v_1\) and \(v_2\) if one of the following conditions holds:
	\begin{enumerate}
		\renewcommand{\theenumi}{\roman{enumi}}
		\renewcommand{\labelenumi}{(\theenumi)}
		\item The e-generalized orthodisk has interior angle \(\pi/2\) at \(v_i\).
		\item The vertex \(v_i\) lies at infinity.
		\item The e-generalized orthodisk has interior angle \(3\pi/2\) at \(v_i\), and the horizontal foliation corresponding to \(\Phi_{\Gamma_j}\) is either parallel or orthogonal to the edges incident at \(v_i\).
	\end{enumerate}
\end{definition}
For \(j = 1, 2, \ldots, 2p-2\), the edge connecting the \(j^{\text{th}}\) and \((j+1)^{\text{th}}\) vertices has endpoints with interior angles \(\pi/2\) and \(3\pi/2\). As in the zigzag case considered by Weber and Wolf \cite{weber1998minimal}, the horizontal foliation corresponding to \(\Phi_{\Gamma_j}\), which realizes the homotopy class \(\Gamma_j\), is aligned with the polygonal edges by construction (see \ref{eqn:Homotopyclassofloop}). Consequently, at a \(3\pi/2\)-vertex the horizontal foliation is either parallel or orthogonal to the incident edges, and condition \((iii)\) of Definition~\ref{defn:admissible_edge} is satisfied. Hence, the edge connecting the \(j^{\text{th}}\) and \((j+1)^{\text{th}}\) vertices is the edge for which the foliation  \(\Phi_{\Gamma_j}\) is admissible.

In conclusion, for any e-conjugate orthodisks with partial symmetry of genus \(p\), there exists at least one pair of edges on which \(\Phi_{\Gamma_j}\) is admissible.  

\subsection{When \(j\in\{-1,\;0,\; 2p-1\}\)} Let $(\mu_0,\, \widetilde\mu_0)\in T_\zeta(T_{\lambda_0,p})$. The support of the Beltrami differentials can be taken to be compactly supported in an arbitrarily small neighborhood of the edge pairs (from the stairs steps).  Therefore for $j\in\{-1,\;0,\; 2p-1\}$, 

\[
\mathrm{Re}\int_{M^{G\eta}_{0,2p+2}} \Phi_{\Gamma_j}\,\mu_0
\;=\;
\mathrm{Re}\int_{M^{G^{-1}\eta}_{0,2p+2}} \Phi_{\Gamma_j}'\,\widetilde\mu_0
\;=\;
0,
\]
and hence for  $j=-1,\,0,\,2p-1,$
\(d\,\mathrm{Ext}_{{G\eta}}(\Gamma_j,\cdot)\big|_{\zeta} = 0\) and \(
d\,\mathrm{Ext}_{{G^{-1}\eta}}(\Gamma_j,\cdot)\big|_{\zeta} = 0.
\)
Thus the maps \(\zeta\mapsto\mathrm{Ext}_{G\eta}(\Gamma_j;\, \zeta)\) and  \(\zeta\mapsto\mathrm{Ext}_{G^{-1}\eta}(\Gamma_j;\, \zeta)\) are constant.  Further, these are 
continuous at \(\lambda_0\).   Since $\lambda_0$ corresponds to the genus 1 reflexive pair, hence 
\[
\mathrm{Ext}_{G\eta}(\Gamma_j;\,\lambda_0)
=
\mathrm{Ext}_{G^{-1}\eta}(\Gamma_j;\,\lambda_0),
\quad \text{for } j=-1,\,0,\,2p-1.
\]

Therefore we have the following: 
\begin{lemma}\label{lemma:existence_open_equal_extremal}
For  every \(\zeta\in T_{\lambda_0,p}\) and $
j=-1,\,0,\,2p-1$,
\(
\mathrm{Ext}_{{G\eta}}(\Gamma_j;\zeta)
\;=\;
\mathrm{Ext}_{{G^{-1}\eta}}(\Gamma_j;\zeta).
\)
\end{lemma}

\subsection{A real analytic submanifold $\mathcal Y$ }
 This subsection establishes the existence of a real analytic submanifold $\mathcal{Y}$ of dimension 1 in $T_{\lambda_0,p}$, which will be used in the next Subsection~\ref{subsec:Finalstep}.

For $p\geq 3$:  on the boundary of $T_{\lambda_0, p}$, there exists a degenerate pair $\zeta_0^p$, corresponding to an e-reflexive pair of partially symmetric generalized orthodisks of genus $p - 1$, denoted by $\zeta^{p-1}$. Let $\zeta_0^{p} = (\zeta^{p-1};\, 0)$ on $\partial\,T_{\lambda_0,p} \times (-\epsilon,\epsilon)$.
Since, $\zeta_0^p$ is lying on the boundary of $T_{\lambda_0,p}$ any small neighborhood of $\zeta_0^p$ say $U_0$ in $\overline{T_{\lambda_0,p}}$ can be identified as $U\times[0,\epsilon)$ for some small $\epsilon>0$ and $U$ is any small neighborhood of $\zeta^{p-1}$ in $T_{\lambda_0,p-1}.$ Moreover, the following is satisfied (by Assumption \ref{assumption:reflexive_pair_genus_p_minus_1}) for $\zeta^{p-1}$:
\[\mathrm{Ext}_{G\eta}(\Gamma_j;\zeta^{p-1})=\mathrm{Ext}_{G^{-1}\eta}(\Gamma_j;\zeta^{p-1}),\quad~j=-1,\,0,\dots,\,p-2,\,2p-3.\]
For $j=p+1,\ p+2,\ldots,\ 2p-4$, at \(\zeta^{p-1}\) this relation is automatically satisfied due to the symmetries in the staircases. We take the map \(\psi:T_{\lambda_0,p}\times(-\epsilon,\epsilon)\rightarrow\R^{p-2}\) defined as
\begin{align*}
	\psi(\zeta,t) := ( 
	&\, \mathrm{Ext}_{G\eta}(\Gamma_1;\zeta) - \mathrm{Ext}_{G^{-1}\eta}(\Gamma_1;\zeta),  \\
	&\, \mathrm{Ext}_{G\eta}(\Gamma_2;\zeta) - \mathrm{Ext}_{G^{-1}\eta}(\Gamma_2;\zeta), \nonumber \\
	&\, \cdots, \nonumber \\
	&\, \mathrm{Ext}_{G\eta}(\Gamma_{p-2};\zeta) - \mathrm{Ext}_{G^{-1}\eta}(\Gamma_{p-2};\zeta))\nonumber
\end{align*}
and form an equation
\begin{equation}\label{equation:analytic_submanifold}
	\psi(\zeta,t)=0,
\end{equation}
and restrict it (after identification) $\psi:U\times(-\epsilon,\epsilon)\rightarrow \R^{p-2} $.  For $\zeta_0^p,$ the following are satisfied:
\begin{align*}
	t_{-1} &= \phi(s_{-1}),\\
	t_{0} &= \phi(s_{0}),\\
	t_j &= \phi(s_j), \quad j \notin \{p-1,\,p,\,p+1\},\\
	t_{p-1} &= t_p = t_{p+1}
	= \phi(s_{p-1}) = \phi(s_p) = \phi(s_{p+1}),
\end{align*}
for some conformal map $\phi$ between the degenerate pair of Riemann spheres
corresponding to $\zeta_0^{p}$. This is exactly the same situation as in \cite[Subsection 5.3, Lemma 5.5]{weber1998minimal}. From the proof of Lemma 5.5 of \cite{weber1998minimal}, the map is differentiable at $\zeta_0^p$,  and if $(\dot\mu,\bar{\dot\mu},v)$ lies in the kernel of $d\psi\vert_{\zeta_0}$, then
\begin{equation*}
	(\dot\mu,\Bar{\dot\mu})=0.
\end{equation*}
In other words, the Jacobian of $\psi$ at $\zeta^p_{0}$ is of full rank. Thus, from the Implicit function theorem there exists an open set in $U\times(-\epsilon,\,\epsilon)$ say $U_0\times(-\epsilon_1, \epsilon_1)$ where the following are satisfied:
\begin{equation}\label{equation:extremal length condition}
	\mathrm{Ext}_{G\eta}(\Gamma_j;(\zeta,t))=\mathrm{Ext}_{G^{-1}\eta}(\Gamma_j;(\zeta,t)),\quad~j=1,\,2,\,\dots,\,p-2.
\end{equation}
Thus, we have a 1-dimensional real analytic manifold $\mathcal{Y}$ such that $\mathcal{Y}$ lies in the zero set of $\psi$. 
\medskip

Note that apart from \eqref{equation:extremal length condition},  from Lemma \ref{lemma:existence_open_equal_extremal}, in $\mathcal{Y}$ the following also holds:
\begin{equation}\label{equation:extremal_length_genus_one_case}
	\mathrm{Ext}_{{G\eta}}(\Gamma_j;\,\zeta)
	\;=\;
	\mathrm{Ext}_{{G^{-1}\eta}}(\Gamma_j;\,\zeta),
	\quad
	j=-1,\,0,\,2p-1
\end{equation}
and in $\mathcal{Y}$,  $H^p_{\Gamma_i}(\zeta)=0$ for all $i\in \{-1,\,0,\,1,\,\dots,\,p-2,\, 2p-1\},\,i\neq p-1$. Thus, similar to \cite[Subsection 5.3]{weber1998minimal}, it follows that $\mathcal{Y}$ acquires the structure of a one-dimensional real analytic submanifold properly embedded in $T_{\lambda_0,p}$. 

For genus $p=2$, note that \eqref{equation:analytic_submanifold} is vacuously solved. Further, from Lemma \ref{lemma:existence_open_equal_extremal} it follows that for every $\zeta\,\in\, T_{\lambda_0,2}$, \eqref{equation:extremal_length_genus_one_case} holds for \(j=-1,\,0,\,3.\)
Therefore, the manifold $T_{\lambda_0,2}$ turns out to be the zero set of real analytic equations; it is a real analytic manifold of dimension $1$.  Take $\mathcal Y\,:=\, T_{\lambda_0,2}$.

\subsection{Final step: finding a reflexive pair of generalized orthodisks}\label{subsec:Finalstep}
This section -- and the article -- concludes by justifying the existence of critical points in $\mathcal{Y}$, and showing that \textit{every critical point of $H^p$ is e-reflexive}.

By Section \ref{sec:height_func},  the function  $H^p$  is proper $C^1$ map on $\mathcal{Y}$. On $\mathcal Y,$
\[H^p(\zeta)=H^p_{\Gamma_{p-1}}(\zeta).\]
For $p\geq 2$, similar to the zigzag, there is at least one  edge  having admissible foliation  \cite[Subsection 6.4]{bardhan2023higher}, \cite{weber1998minimal}, hence there is  \((\dot\mu, \Bar{\dot\mu})\in T_{\zeta} T_{\lambda_0, p}\) for $\zeta\in\,\mathcal{Y}$ such that 
$$\text{sgn}\left(d\mathrm{Ext}_{G\eta}(\Gamma_{p-1};\zeta) (\dot{\mu})  \right) = -\text{sgn}\left( d\mathrm{Ext}_{G^{-1}\eta}(\Gamma_{p-1};\zeta) (\overline{\dot{\mu}})\right).$$
The derivative of $H^p$ is given by: 
\begin{align*}  
D H^{p}\vert_{\zeta}\left(\dot{\mu},\overline{\dot{\mu}}\right) 
&= 2\left(e^{\mathrm{Ext}_{G\eta}(\Gamma_{p-1};\zeta)} - e^{\mathrm{Ext}_{G^{-1}\eta}(\Gamma_{p-1};\zeta)}\right) \\
&\quad \times \left( 
e^{\mathrm{Ext}_{G\eta}(\Gamma_{p-1};\zeta)} 
\left( d\mathrm{Ext}_{G\eta}(\Gamma_{p-1};\zeta) \right)\dot{\mu} 
- e^{\mathrm{Ext}_{G^{-1}\eta}(\Gamma_{p-1};\zeta)} 
\left( d\mathrm{Ext}_{G^{-1}\eta}(\Gamma_{p-1};\zeta) \right)\overline{\dot{\mu}} 
\right) \\
&\quad + 2\left( 
e^{\frac{1}{\mathrm{Ext}_{G\eta}(\Gamma_{p-1};\zeta)}} 
- e^{\frac{1}{\mathrm{Ext}_{G^{-1}\eta}(\Gamma_{p-1};\zeta)}} 
\right) \\
&\quad \times \left( 
- \frac{e^{\frac{1}{\mathrm{Ext}_{G\eta}(\Gamma_{p-1};\zeta)}}}{\mathrm{Ext}^2_{G\eta}(\Gamma_{p-1};\zeta)} 
\left( d\mathrm{Ext}_{G\eta}(\Gamma_{p-1};\zeta) \right)\dot{\mu} 
+ \frac{e^{\frac{1}{\mathrm{Ext}_{G^{-1}\eta}(\Gamma_{p-1};\zeta)}}}{\mathrm{Ext}^2_{G^{-1}\eta}(\Gamma_{p-1};\zeta)} 
\left( d\mathrm{Ext}_{G^{-1}\eta}(\Gamma_{p-1};\zeta) \right)\overline{\dot{\mu}} 
\right).
\end{align*}

If $\zeta$ is not e-reflexive, that implies \(\mathrm{Ext}_{{G\eta}}(\Gamma_{p-1};\zeta)\neq\mathrm{Ext}_{{G^{-1}\eta}}(\Gamma_{p-1};\zeta).\)
Therefore, the both terms of $D H^p\vert_{\zeta}(\dot{\mu},\Bar{\dot{\mu}})$ are strictly positive or strictly negative. Without loss of generality we may assume $DH^p\vert_{\zeta}(\dot{\mu},\Bar{\dot{\mu}})>0$.

 On the other hand, recall that $H^p$ is proper in $\mathcal{Y}$. This implies the  existence of  a point $\Bar{\zeta}$ in $T_{\lambda_0,p}$ such that $DH^p\vert_{\Bar{\zeta}}=0.$   So we have a contradiction.  Therefore, there must exist a point $\zeta$ in $\mathcal{Y}$, such that $\zeta$ is e-reflexive.

Thus, summarizing the discussion, we formalize:
\begin{theorem}
    There exists an e-reflexive pair of generalized orthodisks as in Section \ref{section:formalWeierstrassdataForGenusP} for any genus $p\geq 1$ and hence genus \(p\) Angel surface exists.
\end{theorem}
\bibliographystyle{plainurl}
\bibliography{biblio.bib}

\begin{thebibliography}{1}

\bibitem{bardhan2023higher}
Rivu Bardhan, Indranil Biswas, and Pradip Kumar.
\newblock Higher genus maxfaces with enneper end.
\newblock {\em The Journal of Geometric Analysis}, 34(7), 2024.
\newblock URL: \url{https://doi.org/10.1007/s12220-024-01661-2}.

\bibitem{fujimoriangel2016}
Shoichi Fujimori and Toshihiro Shoda.
\newblock Minimal surfaces with two ends which have the least total absolute curvature.
\newblock {\em Pacific Journal of Mathematics}, 282(1):107--144, 2016.
\newblock \href {https://doi.org/10.2140/pjm.2016.282.107} {\path{doi:10.2140/pjm.2016.282.107}}.

\bibitem{osserman1964annals}
Robert Osserman.
\newblock Global properties of minimal surfaces in {$E^3$} and {$E^n$}.
\newblock {\em Ann. of Math. (2)}, 82(2):340 -- 364, 1964.

\bibitem{sato1996tohoku}
Katsunori Sato.
\newblock Construction of higher genus minimal surfaces with one end and finite total curvature.
\newblock {\em Tohoku Math. J.}, 48(2):229 -- 246, 1996.

\bibitem{schoen1983uniqueness}
Richard~M Schoen.
\newblock Uniqueness, symmetry, and embeddedness of minimal surfaces.
\newblock {\em Journal of Differential Geometry}, 18(4):791--809, 1983.

\bibitem{stein2010complex}
Elias~M Stein and Rami Shakarchi.
\newblock {\em Complex analysis}, volume~2.
\newblock Princeton University Press, 2010.

\bibitem{weberwebsite}
Mathias Weber.
\newblock The angel surfaces, 2018.
\newblock Accessed: 2024-04-22.
\newblock URL: \url{https://theinnerframe.org/2018/06/18/the-angel-surfaces/}.

\bibitem{weber1998minimal}
Matthias Weber and Michael Wolf.
\newblock Minimal surfaces of least total curvature and moduli spaces of plane polygonal arcs.
\newblock {\em Geometric And Functional Analysis GAFA}, 8, 1998.
\newblock \href {https://doi.org/10.1007/s000390050125} {\path{doi:10.1007/s000390050125}}.

\bibitem{weber2002teichmuller}
Matthias Weber and Michael Wolf.
\newblock Teichm{\"u}ller theory and handle addition for minimal surfaces.
\newblock {\em Annals of mathematics}, pages 713--795, 2002.

\end{thebibliography}
\end{document}